%% file: main.tex
\documentclass[11pt,reqno]{amsart}
\usepackage[T1]{fontenc}
\usepackage{lmodern}
\usepackage{amsmath,amssymb,mathtools}
\usepackage{microtype}
\usepackage[a4paper,margin=28mm,headheight=14pt,headsep=7mm]{geometry}
\usepackage{booktabs,array,longtable}
\usepackage{enumitem}
\usepackage{tikz}
\usetikzlibrary{arrows.meta,calc}
\usepackage{hyperref}
\hypersetup{hidelinks,pdftitle={Coxeter Rigidity and ADE Classification of Fibonacci-Generated Tensor Categories},pdfauthor={HUIXUAN HE and ZHENGWEI LIU},pdfsubject={Fibonacci-generated tensor categories, Coxeter rigidity, marked categorification in finite and infinite rank, and finite ADE classification}}
\usepackage{aliascnt}
\usepackage[nameinlink,noabbrev]{cleveref}
\newtheorem{theorem}{Theorem}[section]
\newaliascnt{proposition}{theorem}
\newtheorem{proposition}[proposition]{Proposition}
\aliascntresetthe{proposition}
\newaliascnt{lemma}{theorem}
\newtheorem{lemma}[lemma]{Lemma}
\aliascntresetthe{lemma}
\newaliascnt{corollary}{theorem}
\newtheorem{corollary}[corollary]{Corollary}
\aliascntresetthe{corollary}
\theoremstyle{definition}
\newaliascnt{definition}{theorem}
\newtheorem{definition}[definition]{Definition}
\aliascntresetthe{definition}
\theoremstyle{remark}
\newaliascnt{remark}{theorem}
\newtheorem{remark}[remark]{Remark}
\aliascntresetthe{remark}
\crefname{theorem}{Theorem}{Theorems}
\crefname{proposition}{Proposition}{Propositions}
\crefname{lemma}{Lemma}{Lemmas}
\crefname{corollary}{Corollary}{Corollaries}
\crefname{definition}{Definition}{Definitions}
\crefname{remark}{Remark}{Remarks}
\numberwithin{equation}{section}
\newcommand{\Fib}{\mathrm{Fib}}
\newcommand{\KFib}{K_{\mathrm{Fib}}}
\newcommand{\TTthree}{\mathcal{TT}_3}
\input{companion_refs}
\newcommand{\cC}{\mathcal C}
\newcommand{\cD}{\mathcal D}
\newcommand{\cE}{\mathcal E}
\newcommand{\cM}{\mathcal M}
\newcommand{\eobj}[1]{\mathbf e_{#1}}
\newcommand{\Czero}{\mathcal C_0}
\newcommand{\HH}{H_{\mathrm{Fib}}}
\newcommand{\one}{\mathbf1}
\newcommand{\id}{\operatorname{id}}
\newcommand{\Id}{\operatorname{Id}}
\newcommand{\Hom}{\operatorname{Hom}}
\newcommand{\End}{\operatorname{End}}
\newcommand{\im}{\operatorname{im}}
\newcommand{\rank}{\operatorname{rank}}
\newcommand{\FPdim}{\operatorname{FPdim}}
\renewcommand{\Vec}{\mathsf{Vec}}

\newcommand{\TT}{\mathbb T}
\newcommand{\ZZ}{\mathbb Z}
\newcommand{\QQ}{\mathbb Q}
\newcommand{\CC}{\mathbb C}
\newcommand{\Cx}{\mathbb C^{\times}}
\newcommand{\RR}{\mathbb R}
\newcommand{\mix}{\mathrm{mix}}

\newcommand{\eps}{\varepsilon}
\newcommand{\ph}{\varphi}
\newcommand{\cell}[3]{e(#1,#2,#3)}
\newcommand{\arxiv}[1]{\href{https://arxiv.org/abs/#1}{arXiv:#1}}
\newcommand{\doi}[1]{\href{https://doi.org/#1}{doi:#1}}
\setlist[enumerate]{label=\textup{(\roman*)},leftmargin=2.3em,itemsep=2pt,topsep=4pt}
\allowdisplaybreaks[1]

\newcommand{\ord}{\operatorname{ord}}

\newcommand{\iso}{\xrightarrow{\ \sim\ }}
\newaliascnt{example}{theorem}

\aliascntresetthe{example}
\crefname{example}{Example}{Examples}
\title[Coxeter rigidity and ADE classification]{Coxeter Rigidity and ADE Classification\\ of Fibonacci-Generated Tensor Categories}
\author{HUIXUAN HE}
\address{H. He, Department of Mathematics,
Tsinghua University, Beijing, China}
\email{hhx22@mails.tsinghua.edu.cn}
\author{ZHENGWEI LIU}
\address{Z. Liu, Yau Mathematical Sciences Center
and Department of Mathematics, Tsinghua University,
and Beijing Institute of Mathematical Sciences and Applications,
Beijing, China}
\email{liuzhengwei@mail.tsinghua.edu.cn}
\date{October 23, 2026}
\subjclass[2020]{Primary 18M20; Secondary 20F55, 46L37}
\keywords{Fibonacci category, fusion category, Coxeter group, Hecke algebra, module functor, unitarizability, Davis complex, gauge equivalence}
\begin{document}
\begin{abstract}
We propose a program to classify tensor categories generated by
prescribed seed tensor categories. We carry out this program for
finitely many marked copies of the positive Fibonacci category with
pairwise nonisomorphic generating simple objects, in the setting of
split locally finite semisimple rigid $\CC$-linear tensor categories
with simple unit. We show that the Fibonacci fusion relation forces
the fusion rules to be those of a specialized Coxeter--Hecke ring.
The systems admitting categorification have Coxeter exponents
$2,3,\infty$. For every such system we construct a tensor category
and prove uniqueness up to marked tensor equivalence. Its simple
objects are indexed by the
Coxeter group, which may be infinite. Finite rank occurs precisely
for disjoint unions of finite ADE diagrams, giving the ADE
classification of Fibonacci-generated fusion categories. These
categories are unitarizable, with no unitary structure assumed on
the candidate. In finite ADE type, allowing Yang--Lee seeds gives
one independent tensor-branch choice on each connected component.
\end{abstract}
\maketitle

\section{Introduction}\label{sec:introduction}
We propose a program to classify composites of prescribed tensor
categories. Given a finite family of seed tensor categories, we seek
tensor categories generated by fully faithful tensor images of these
seeds, with the embeddings retained as \emph{markings}. The seeds may
be copies of the same category or may be different. Classification is
up to tensor equivalence compatible with the markings. The problem is
to determine which mixed tensor product relations can occur, which
fusion rules admit a compatible tensor structure, and how much freedom
remains in the associativity constraints once the seed tensor
structures are fixed. A further question is which composites have
finitely many simple objects.

In this paper we carry out the program for positive Fibonacci seeds.
The Fibonacci category $\Fib$ has two simple objects, $\one$ and $\tau$,
with
\[
 \tau\otimes\tau\cong\one\oplus\tau.
\]
We call a tensor category generated by copies of $\Fib$
\emph{Fibonacci-generated}, or a \emph{composite} of Fibonacci
categories. Generation is under tensor products, finite direct sums,
and summands. We consider finitely many copies with pairwise
nonisomorphic generating simples. Each copy carries the fixed
positive Fibonacci tensor structure, including its associator.
The ambient categories are split locally finite semisimple rigid
$\CC$-linear tensor categories with simple unit; they may have
infinitely many simple isomorphism classes.

The composition problem is closely related to subfactor composition.
For a chain $N\subset P\subset M$, the two components give bimodule
categories on the common $P$ side, and their interaction records the
relative position of the inclusions. Bisch and Jones developed the
Fuss--Catalan framework for free compositions and intermediate
subfactors \cite{BischJones}; Bisch and Haagerup analyzed group-type
compositions \cite{BischHaagerup} and posed the problem of which fish
graphs arise from an $A_3$--$A_4$ composition. Here $A_3$ and $A_4$
denote the principal graphs of the component subfactors, whose even
parts are $\Vec_{\mathbb Z/2\mathbb Z}$ and $\Fib$, respectively.
The fish-graph problem, posed in 1994, was resolved by Liu, who also
classified the standard invariants of $A_4$--$A_4$ compositions
\cite{Liu}. Izumi--Morrison--Penneys formulated this problem in terms
of categorical quotients and developed a complementary approach through
jellyfish relations. In their unitary setting, the
$A_4$--$A_4$ composition problem translates into the problem of composing
two copies of $\Fib$, or equivalently, studying the corresponding tensor
quotients of $\Fib*\Fib$: a composite subfactor gives such an even
part, and algebra objects in these quotients reconstruct the composite
subfactors \cite[expanded preprint, Section~5.1]{IMP}.

Two basic composites of $\Fib$ with itself are the Deligne product
and the free product.
In the former the factors commute; in the latter, alternating words
remain simple. Between them lies the six-simple-object category
$\TTthree$, the dual even part of the $2D2$ subfactor
\cite{IMP,Liu,EM}. Its two Fibonacci generators satisfy
$aba\cong bab$. This basic finite noncommuting example
suggests that mixed relations should be organized by Coxeter systems.

We first recover a Coxeter system from the fusion rules and then
construct and classify the corresponding tensor categories.
For distinct marked positive Fibonacci generators indexed by a finite set $S$,
the possible composites are parametrized by Coxeter matrices with
off-diagonal entries $2,3,\infty$. Each matrix determines one
complex-linear tensor category up to marked tensor equivalence.
Its simple objects are indexed by the Coxeter group $W$, and its
fusion rules are the standard-basis Hecke rules with quadratic
relation $h_s^2=1+h_s$. The group $W$ may be infinite, even though the
generating set is finite. The category is a fusion category precisely
when the diagram is a disjoint union of finite ADE diagrams.

The diagram records the generating subcategories and their pairwise
relations: each vertex represents a copy of $\Fib$, and the edges and
their labels encode the Coxeter exponents. Here \emph{Coxeter rank}
is $|S|$, the number of marked generators, whereas \emph{categorical
rank} is $|W|$, the number of simple isomorphism classes. A single
vertex gives $\Fib$, the edge $A_2$ gives
$\TTthree$, and the three-vertex path $A_3$ gives a category with
$24$ simple objects. Already in the last example, a complete table of
$F$-symbols contains more than a million entries
(\cref{subsec:A3}).

\subsection{The main theorems}
The argument proceeds from based-ring rigidity to existence and
marked uniqueness of tensor categorifications. Together these give
\cref{cor:classification} for finite or infinite Coxeter groups;
the classical finiteness criterion then yields the ADE case in
\cref{cor:ADE}.
Fix the positive complex-linear Fibonacci model, including its
associator, unit, duality maps, and fusion decompositions. The
generating object satisfies
\begin{equation}\label{eq:fib-rule}
\tau\otimes\tau\cong\one\oplus\tau.
\end{equation}
\medskip\noindent\textbf{I. Coxeter rigidity and the categorical restriction.}
Write $\HH(W)$ for the standard-basis Coxeter--Hecke ring with
generator multiplication \eqref{eq:hecke-rule}.

\begin{theorem}[Coxeter rigidity]\label{thm:intro-rigidity}
Let $R$ be a rigid based ring generated as a based ring by distinct basis elements $x_s$, $s\in S$, where $S$ is finite and $x_s^2=1+x_s$. Then $R\cong\HH(W)$ for a Coxeter system $(W,S)$, by an isomorphism taking $x_s$ to $h_s$. Its matrix is determined by the marked two-generator subrings. At the level of based rings, the off-diagonal entries may take arbitrary values in $\{2,3,4,\ldots\}\cup\{\infty\}$.
\end{theorem}

The proof is given in \cref{thm:ring-rigidity} and applies to based
rings of arbitrary rank. It recovers the Coxeter system from the
based-ring axioms alone, before any associator is chosen.
For tensor categories, the arithmetic
obstruction of \cite[Theorem~\PointedArithmeticTheorem]{FibPointed}
restricts the possible exponents.
Here and below ``locally finite semisimple'' means that Hom spaces are finite-dimensional and every object is a finite direct sum of simple objects; the set of simple isomorphism classes may be infinite. We use \emph{fusion rules} also in this infinite-rank setting, with finite support in each tensor product.

\begin{corollary}[Categorifiability restriction]\label{cor:intro-exponents}
Let $\cD$ be a split locally finite semisimple rigid $\CC$-linear tensor
category with simple unit, generated under tensor products, finite direct
sums, and summands by finitely many pairwise nonisomorphic simple objects
$Y_s$ satisfying $Y_s^{\otimes2}\cong\one\oplus Y_s$.
Then $K_0(\cD)\cong\HH(W)$ for a marked Coxeter system $(W,S)$, with
$[Y_s]\mapsto h_s$, and every off-diagonal exponent belongs to
$\{2,3,\infty\}$.
\end{corollary}

The proof in \cref{cor:exponent-restriction} applies the arithmetic
obstruction to each finite dihedral tensor subcategory.
The construction and uniqueness results below fix the positive
Fibonacci tensor branch on every generator.

\medskip\noindent\textbf{II. Algebraic construction.}
Let $(W,S)$ satisfy
\begin{equation}\label{eq:exponents}
|S|<\infty,\qquad m_{st}\in\{2,3,\infty\}\quad(s\ne t).
\end{equation}
Write $\ell$ for Coxeter length. The based ring $\HH(W)$ has basis $\{h_w:w\in W\}$ and generator multiplication
\begin{equation}\label{eq:hecke-rule}
h_sh_w=\begin{cases}
h_{sw},&\ell(sw)=\ell(w)+1,\\
h_{sw}+h_w,&\ell(sw)=\ell(w)-1.
\end{cases}
\end{equation}
It is the Hecke specialization in the standard basis, with quadratic relation $h_s^2=1+h_s$.

\begin{theorem}[Algebraic construction]\label{thm:construction}
For every system satisfying \eqref{eq:exponents}, there exists a split locally finite semisimple rigid $\CC$-linear tensor category $\Czero(W)$ with simple unit and the following properties.
\begin{enumerate}
\item Its simple objects are $X_w$, $w\in W$, with $X_e=\one$ and $X_w^{\vee}\cong X_{w^{-1}}$.
\item Multiplication by $X_s$ satisfies the direct-sum form of \eqref{eq:hecke-rule}; hence $K_0(\Czero(W))\cong\HH(W)$ as based rings.
\item For each $s\in S$ there is a fully faithful tensor functor $\iota_s^0:\Fib\to\Czero(W)$ taking $\tau$ to $X_s$, and these images generate $\Czero(W)$.
\item The category has finitely many simple objects if and only if $W$ is finite. In that case
\begin{equation}\label{eq:dimension-main}
\operatorname{rank}\Czero(W)=|W|,\qquad
\FPdim(X_w)=\ph^{\ell(w)},\qquad
\FPdim\Czero(W)=\sum_{w\in W}\ph^{2\ell(w)},
\end{equation}
where $\ph=(1+\sqrt5)/2$.
\end{enumerate}
\end{theorem}

The construction is completed in \cref{sec:construction}; it realizes
every Coxeter system allowed by \cref{cor:intro-exponents}, including
infinite groups. Proposition~\ref{prop:unitary-realization} then gives
a unitary realization. It remains to determine whether these marked
fusion rules admit any other tensor structure.

\medskip\noindent\textbf{III. Marked algebraic uniqueness.}
\begin{definition}[Marked categorification]\label{def:marked}
A marked categorification of $\HH(W)$ by the fixed $\Fib$ consists of
a split locally finite semisimple rigid $\CC$-linear tensor category
$\cD$ with simple unit, a based-ring isomorphism
$\kappa:K_0(\cD)\iso\HH(W)$, and fully faithful tensor functors
$\iota_s:\Fib\to\cD$ taking $\tau$ to the simple object $Y_s$ with
$\kappa([Y_s])=h_s$.
A marked equivalence from $(\cD,\kappa,(\iota_s)_{s\in S})$ to
$(\cD',\kappa',(\iota_s')_{s\in S})$ consists of a complex-linear
tensor equivalence
\[
F:\cD\longrightarrow\cD',\qquad \kappa'\circ K_0(F)=\kappa,
\]
together with monoidal natural isomorphisms
\[
\eta_s:F\circ\iota_s\overset{\sim}{\Longrightarrow}\iota_s'
\qquad(s\in S).
\]
\end{definition}

\begin{theorem}[Fixed-seed algebraic uniqueness]\label{thm:main}
Every marked categorification in \cref{def:marked} is marked complex-linearly tensor equivalent to $\Czero(W)$. No dagger, pivotal, or spherical structure is assumed on the candidate. The group $W$ need not be finite.
\end{theorem}

\medskip\noindent\textbf{IV. Classification and finite ADE type.}
The preceding results classify the marked tensor categories. Their
finite part is determined by the classical finite Coxeter diagrams.

\begin{corollary}[Finite ADE composites]\label{cor:ADE}
For each disjoint union of finite ADE diagrams, the ring $\HH(W)$ has exactly one marked complex-linear categorification by the fixed $\Fib$, up to marked tensor equivalence. It is unitarizable and has the dimensions in \eqref{eq:dimension-main}. Conversely, the construction of \cref{thm:construction} has finite rank only for these diagrams.
\end{corollary}
\begin{proof}
By the finiteness criterion in \cref{subsec:geometric-representation},
a Coxeter system satisfying \eqref{eq:exponents} is finite exactly when
its diagram is a disjoint union of finite ADE diagrams. Existence and dimensions follow from \cref{thm:construction},
marked uniqueness from \cref{thm:main}, and unitarizability from
\cref{prop:unitary-realization}.
\end{proof}

\begin{corollary}[Classification with marked positive generators]
\label{cor:classification}
Fix a finite label set $S$. Consider split locally finite semisimple rigid
$\CC$-linear tensor categories with simple unit and specified fully
faithful tensor functors $\iota_s:\Fib\to\cD$, whose images tensor-generate
$\cD$ and whose objects $\iota_s(\tau)$ are pairwise nonisomorphic.
Up to complex-linear tensor equivalence equipped with monoidal natural
isomorphisms on all labelled seed functors, these categories are classified
by Coxeter matrices on $S$ with off-diagonal entries in $\{2,3,\infty\}$.
The representative is $\Czero(W)$. It is unitarizable and is a fusion
category exactly when the Coxeter diagram is a disjoint union of finite
ADE diagrams.
\end{corollary}
\begin{proof}
Based-ring rigidity identifies $K_0(\cD)$ with $\HH(W)$, with the
marked generators corresponding to the standard generators. The marked
Coxeter matrix is then read off from $m_{st}=\ord_W(st)$, and
\cref{cor:exponent-restriction} gives the allowed exponents.
Theorem~\ref{thm:construction} supplies every representative, while
\cref{thm:main} supplies its marked uniqueness. An equivalence preserving
the seed functors preserves their two-generator based subrings and hence
the matrix. Unitarizability and the finite-type assertion follow from
\cref{prop:unitary-realization,cor:ADE}.
\end{proof}

The uniqueness assertion fixes the local \emph{tensor} branch. The
two-simple-object Fibonacci based ring also admits the Yang--Lee branch
\cite[Section~3]{BookerDavydov}. In finite ADE type, allowing both
branches amounts to choosing one branch on each connected component;
\cref{cor:ADE-branches} gives the resulting marked classification.
The positive-seed classification above applies also to infinite $W$.

\subsection{Proof strategy}
The fusion rule first determines a group. Multiplication by a
Fibonacci basis element pairs the basis into two-point orbits, with
one additional diagonal term on each pair. Associativity makes the
left and right pairings commute. Their regular permutation actions
recover a group, and the diagonal terms recover its descent sets.
The exchange condition then identifies it as a Coxeter group. The
arithmetic theorem of the companion paper excludes finite exponents
other than two and three.

A general regular realization guides the construction. In the split
locally finite semisimple setting, fusion rules determine a linear
category and its right multiplication functors. Endofunctors equipped
with natural commutation isomorphisms with these functors form a monoidal
category; its morphisms preserve these isomorphisms. This ambient category
requires no choice of associator. Every categorification embeds fully
faithfully by left tensoring, with its associator supplying the
commutation isomorphisms. Finding tensor subcategories with the prescribed fusion
rules is therefore a universal approach to categorification. Here we
can impose further compatibility: for each colour, the auxiliary
category decomposes into blocks with two simples, each carrying a copy
of the regular right $\Fib$-module. The fixed Fibonacci tensor data thus supply
complete right module structures. Requiring the commutation isomorphisms
to satisfy the module-functor axioms makes them module-functor constraints. This regular decomposition is special to the seed;
other seeds can be treated similarly when suitable module data are
available. The fixed Fibonacci module structures thus define the ambient
category in which we construct the composite tensor category.

We implement this approach on $\cM=\Vec_W^{\mathrm{fin}}$, the linear
category of finite-support $W$-graded vector spaces, with complete
left and right Fibonacci actions.
We build the two finite braid intertwiners and use parabolic coordinates
to verify their compatibility with every right colour. Inside the
resulting category of module endofunctors, evaluation at the identity
separates reduced words. The braid isomorphisms, together with the
splittings $L_s^2\cong\Id\oplus L_s$, express every finite word in
the generating functors as a finite direct sum of objects represented
by reduced words. These facts establish semisimplicity. Rigidity
follows from the duality maps of the seed actions and their extension
to finite direct sums and direct summands. The same formulas admit
a positive Hilbert-space realization.

For uniqueness, we realize a candidate by its regular actions on
the same $\cM$ and choose identifications of its one-colour regular
blocks with the fixed models. In these coordinates, mixed interchanges
and normalized braid morphisms are described by nonzero face parameters
$h,k,\theta$. Changing the block identifications rescales the coloured
edges. Since an edge can belong to several faces, these rescalings
must be chosen simultaneously. The space
\[
 Y=(\Sigma\times\Sigma)/W,
\]
where $\Sigma$ is the Davis complex introduced in
\cref{subsec:davis-prelim} and $W$ acts diagonally, organizes this compatibility. We prove that its quotient cells have
ball characteristic domains and that $Y$ is contractible. Two-generator
marked equivalences give the scalar equations on mixed three-cells; integral
boundary elimination makes them sufficient for every two-cycle.
A spanning tree then yields simultaneous edge rescalings.
The fully faithful regular realizations have the same generated tensor
image; lifting the tensor structure from this image recovers the tensorator and the
marked tensor equivalence. Working over $\ZZ$ preserves finite-order scalars as well
as continuous parameters.

\subsection{Relation to earlier work}
The companion paper \cite{FibPointed} studies Fibonacci--pointed
composites. Its formal-codegree exclusion is our sole proof input from
that paper; its geometric Coxeter rigidity and graphical two-generator
uniqueness give independent proofs by different methods.

Yamagami constructs free products through recursively defined triangular
morphism spaces and associativity maps \cite[Sections~1--3]{Yamagami}.
Valvekens uses compatible left and right actions on graded Hilbert spaces
inside a category of module endofunctors \cite[p.~2]{Valvekens}, the direct
precedent for our regular-action approach. We incorporate finite Coxeter
braid intertwiners, algebraic reconstruction, and global marked equivalence.

Section~\ref{sec:structure} identifies $\Czero(A_2)$ with $\TTthree$.
The commuting two-factor case also follows from the twisted Deligne
product theorem \cite[Theorem~1.2(i)]{ENOtwisted}, since $\Fib$ has trivial
universal grading.\footnote{If $\tau$ has degree $g$, its two fusion
summands give $g^2=e=g$.}
Our explicit local normalization provides coordinates for the global
uniqueness proof.

Our $A_3$ example contradicts \cite[Lemma~4.7, p.~819, and
Sublemma~4.7.4, pp.~830--831]{EM} for rank-preserving quotients, meaning
quotients with pairwise nonisomorphic marked generators
\cite[Definition~2.6]{EM}. Propositions~\ref{prop:free-product-quotient}
and~\ref{prop:A3-comparison} give the details. These intermediate claims enter the proofs of
\cite[Lemma~4.8 and Theorem~4.1]{EM}, so the corresponding steps require
additional justification. The cyclic classification assumes a tensor
autoequivalence cyclically permuting the generators, a symmetry absent
from our $A_3$ example.

The elementary actions relate to the nonnegative integer matrix
representation calculations of \cite{BookerDavydov,ADC}; the Coxeter
and Davis-complex tools come from \cite{Davis}. In comparison with
Soergel categorification \cite{EW}, our categories are semisimple and
rigid, with simples realizing the specialized standard Hecke basis.
Reutter's theorem \cite{Reutter} compares unitary realizations of
unitarizable fusion categories; our algebraic uniqueness establishes
unitarizability of the marked candidates.

\subsection{Organization}
Section~\ref{sec:preliminaries} fixes the Fibonacci tensor data and
introduces Coxeter combinatorics, the reflection representation, the
Davis complex, and the specialized Hecke algebra. Section~\ref{sec:rigidity} proves based-ring rigidity and
the categorifiability restriction. Section~\ref{sec:A2} constructs the
six-object $A_2$ model. Sections~\ref{sec:braids}--\ref{sec:construction}
extend the actions and braid maps to arbitrary Coxeter systems,
reconstruct the categories, and establish unitarizability.
Section~\ref{sec:structure} gives structural consequences and examples,
including the $A_3$ fusion data, direct associator calculation, and
comparison with three-generator quotients in \cref{subsec:A3}.
Sections~\ref{sec:regular}--\ref{sec:rank-two} reduce the tensor-equivalence
problem to a comparison of regular actions and classify the local
bimodule data. Sections~\ref{sec:geometry}--\ref{sec:global}
encode the local parameters and edge rescalings in cellular chains, construct a simultaneous
normalization, and recover the tensorator of the marked equivalence.
Section~\ref{sec:consequences} brings these results together in the
finite ADE classification, then treats unitary structures, tensor
branches, infinite type, and the matrix of edge rescalings.
Section~\ref{sec:discussion} compares the role of ADE here with earlier
classifications and discusses composites of other seed categories.
The appendices give the finite matrix certificates for the local braid
construction and two-generator uniqueness.

\section{Fibonacci categories and Coxeter geometry}\label{sec:preliminaries}
We fix the Fibonacci tensor data and the Coxeter background used
throughout the paper. The reflection representation explains finite
type and supplies the cells of the Davis complex; the specialized
Hecke algebra gives the fusion rules. All later local matrices use
the channel order fixed below.

\subsection{The fixed Fibonacci tensor structure}
All categories and functors are complex-linear. A locally finite semisimple category in this paper has finite-dimensional Hom spaces and finite direct-sum decompositions; it may have infinitely many simple isomorphism classes. Tensor categories have simple unit and both left and right duals. We write $X^{\vee}$ for an object dual and reserve $f^{\dagger}$ for an adjoint in a specified unitary realization. Neither a pivotal structure nor an identification $X\cong X^{\vee\vee}$ is part of the candidate data. Composition is read from right to left.

Recall $\ph=(1+\sqrt5)/2$ and write $\KFib=\QQ(\sqrt5)$. Fix a skeletal model of $\Fib$ with strict unit and specified algebraic fusion decompositions. Set
\begin{equation}\label{eq:U}
a=\ph^{-1},\qquad b=\ph^{-1/2},\qquad
U=F^{\tau\tau\tau}_{\tau}=
\begin{pmatrix}a&b\\b&-a\end{pmatrix},
\end{equation}
where the intermediate channels are ordered as $(\one,\tau)$. All scalar associator blocks are one. The relations
\begin{equation}\label{eq:seed-identities}
a=b^2,\qquad a^2+b^2=1,\qquad U^{-1}=U
\end{equation}
will be used algebraically. The full seed includes its unit and duality maps, as well as both inclusions and projections in $\tau^2\cong\one\oplus\tau$. These algebraic data determine the regular actions used below. Their positive realization is established in \cref{prop:unitary-realization}. For the algebraic classification of tensor categories with two simple objects, see \cite[Section~3]{BookerDavydov}.

If $i_{\one}:\one\to\tau^2$ and $p_{\one}:\tau^2\to\one$ are the specified unit-channel inclusion and projection, one may choose the evaluation and coevaluation of $\tau$ to be $\sqrt\ph\,p_{\one}$ and $\sqrt\ph\,i_{\one}$. The two snake identities reduce to $\ph a=1$. The associator model for this two-simple-object category is the standard algebraic one of \cite[Section~3]{BookerDavydov}.

For a semisimple decomposition, inclusions $i_j$ and projections $p_j$ satisfy $p_ji_k=\delta_{jk}\id$ and $\sum_j i_jp_j=\id$. We use these algebraic identities, together with $U^{-1}=U$, in the comparison of tensor structures.

\subsection{Coxeter systems, length and parabolic subgroups}\label{subsec:coxeter-prelim}
Throughout the Coxeter preliminaries, $S$ is finite and the
off-diagonal exponents may be arbitrary elements of
$\{2,3,4,\ldots,\infty\}$. The restriction \eqref{eq:exponents}
is imposed for the categorical construction and uniqueness.
A Coxeter system has presentation
\[
 W=\langle S\mid s^2=e,\ (st)^{m_{st}}=e\text{ when }m_{st}<\infty\rangle.
\]
Its diagram has vertex set $S$, with an edge between $s$ and $t$
when $m_{st}\ge3$; an edge of exponent $3$ is unlabelled, and all
other edges carry their exponent, including $\infty$.
The length $\ell(w)$ is the minimum length of a word in $S$ representing $w$; such a word is \emph{reduced}. Multiplication by one generator changes length by exactly one. Indeed, the triangle inequality bounds the change by one, and the homomorphism sending every generator to $-1$ fixes the parity of word length. A generator $s$ is a left descent of $w$ when $\ell(sw)<\ell(w)$; right descents are defined similarly. Put
\begin{equation}\label{eq:descents}
 \eps_s(w)=\mathbf1_{\ell(sw)<\ell(w)},\qquad
 \eta_r(w)=\mathbf1_{\ell(wr)<\ell(w)}.
\end{equation}

We use the exchange condition: if $t_1\cdots t_k$ is reduced and $s$ is a left descent of its product, then $s(t_1\cdots t_k)$ is represented by the word obtained by deleting one $t_i$. Conversely, if a group is generated by distinct nonidentity involutions
and satisfies this exchange condition for every shortest expression,
then it is a Coxeter system, with $m_{st}=\ord(st)$
\cite[Theorem~3.3.4]{Davis}, for finite or infinite groups.
We also use the reduced-word theorem: any two reduced expressions are related by braid moves. These statements and the parabolic facts below are standard; see \cite[Chapters 3--4]{Davis}.

For $J\subseteq S$, the subgroup $W_J=\langle J\rangle$ is a standard parabolic subgroup. It is \emph{spherical} if it is finite. Each coset $W_Jd$ has a unique shortest representative $d$, characterized by the absence of left descents in $J$. Every element of that coset is uniquely $ud$, with $u\in W_J$, and
\begin{equation}\label{eq:length-add}
 \ell(ud)=\ell(u)+\ell(d).
\end{equation}
Double cosets $W_JdW_K$ also have unique shortest representatives, characterized by the absence of left $J$-descents and right $K$-descents. In the double-coset arguments below, we apply \eqref{eq:length-add} to each constituent coset $W_Jx$, taking $x$ to be the shortest representative of that coset.

\subsection{The geometric representation}\label{subsec:geometric-representation}
Coxeter generators admit a faithful realization as linear reflections.
Let $V=\bigoplus_{s\in S}\RR\alpha_s$ and equip it with the symmetric
bilinear form
\begin{equation}\label{eq:coxeter-form}
 B(\alpha_s,\alpha_s)=2,\qquad
 B(\alpha_s,\alpha_t)=-2\cos\!\left(\frac{\pi}{m_{st}}\right)
 \quad(s\ne t),
\end{equation}
where $\cos(\pi/\infty)=1$. The action of a generator is
\begin{equation}\label{eq:coxeter-reflection}
 s(v)=v-B(v,\alpha_s)\alpha_s.
\end{equation}
It fixes the hyperplane $B(v,\alpha_s)=0$ pointwise, sends
$\alpha_s$ to $-\alpha_s$, and preserves $B$. These reflections
satisfy the Coxeter relations and give a faithful action of $W$ on $V$
\cite[Section~6.12, especially Corollary~6.12.11]{Davis}.
This is the standard geometric representation, called the
\emph{canonical representation} in that reference.
For $m_{st}<\infty$, the form is positive definite on
$\RR\alpha_s+\RR\alpha_t$; the two reflecting lines bound a sector
of angle $\pi/m_{st}$, and their product is a rotation of order $m_{st}$.

The form $B$ is positive definite exactly when $W$ is finite
\cite[Theorem~6.12.9]{Davis}; in general it may be degenerate or
indefinite. For a spherical subset $J$, its restriction to
$V_J=\bigoplus_{j\in J}\RR\alpha_j$ is positive definite, and $W_J$
acts there as a finite Euclidean reflection group with chamber
\[
 C_J=\{v\in V_J:B(v,\alpha_j)\ge0\text{ for all }j\in J\}.
\]
These finite reflection groups supply the cells introduced below.

Under \eqref{eq:exponents}, an exponent-$\infty$ pair gives the
singular principal block $\left(\begin{smallmatrix}2&-2\\-2&2\end{smallmatrix}\right)$,
so positive definiteness excludes such pairs. When all exponents are
$2$ or $3$, the Gram matrix is $2I-A_\Gamma$, where $A_\Gamma$ is
the adjacency matrix of the Coxeter diagram. Thus $W$ is finite if
and only if $\rho(A_\Gamma)<2$, and the connected finite diagrams are
$A_n,D_n,E_6,E_7,E_8$ \cite[Appendix~C]{Davis}.
Consequently the finite systems satisfying \eqref{eq:exponents} are
precisely those whose diagrams are disjoint unions of finite ADE diagrams.

\subsection{The Davis complex}\label{subsec:davis-prelim}
Throughout, $\Sigma$ denotes the \emph{Davis complex} with its Coxeter cell
structure. Its barycentric subdivision is the order complex of the
spherical cosets \cite[Sections~7.2--7.3]{Davis}.

Let
$\mathcal S=\{J\subseteq S:W_J\text{ is finite}\}$, including
$\varnothing$. For $J\in\mathcal S$, choose $p_J\in C_J$ with
$B(p_J,\alpha_j)=1$ for every $j\in J$, and let $P_J$ be the
convex hull of $W_Jp_J$ in $V_J$. This polytope has dimension $|J|$;
for $J=\varnothing$ it is a point.

For each coset $C=wW_J$ take a copy $wP_J$ of $P_J$, labelling the
vertex corresponding to $up_J$ by $wu$ for $u\in W_J$. The nonempty
faces of this copy have as their vertex sets precisely the spherical
cosets contained in $C$. Glue faces with the same vertex set by
matching their vertex labels.
Left multiplication on the labels gives the $W$-action on $\Sigma$.
Its vertices are the elements of $W$, and its one-skeleton is the
coloured Cayley graph, with edges $w\longleftrightarrow ws$
\cite[Proposition~7.3.4]{Davis}.
A finite exponent $m_{st}$ gives a $2m_{st}$-gon, a square for
$m_{st}=2$ and a hexagon for $m_{st}=3$. An exponent-$\infty$ pair
contributes no two-cell of this type. In \cref{sec:geometry} we use
products of these cells to construct the diagonal quotient.

The barycentric subdivision of $\Sigma$ gives a convenient description
of its fundamental chamber. The face poset of $\Sigma$ is the set
$\mathcal P$ of spherical cosets ordered by inclusion. Its order
complex has a simplex for each finite nonempty chain of cosets.
Placing each coset at the barycentre of its cell (the average of its
vertices) identifies that simplex with the simplex joining the
barycentres of the corresponding nested faces. For example, in
Figure~\ref{fig:davis-A2} the chain
$\{e\}\subset W_{\{s\}}\subset W_{\{s,t\}}$ gives the triangle
joining $e$, the midpoint of the edge $[e,s]$, and the centre of the
hexagon.

The chains of standard cosets $W_J$ form a subcomplex of this
subdivision, the \emph{Davis chamber} $\mathcal K=|\mathcal S|$.
Here $|\mathcal S|$ denotes the realization of the order complex of
the spherical subsets $J$, ordered by inclusion. Every chain of
spherical cosets is a left translate of a unique chain of standard
cosets, so $\Sigma/W\cong\mathcal K$
\cite[Lemma~7.2.3 and Theorem~7.2.4]{Davis}.
The least element $\varnothing$ makes $\mathcal K$ a cone, hence
contractible.

\begin{figure}[ht]
\centering
\begin{tikzpicture}[scale=.85,every node/.style={font=\small}]
 \coordinate (v0) at (0:1.65);
 \coordinate (v1) at (60:1.65);
 \coordinate (v2) at (120:1.65);
 \coordinate (v3) at (180:1.65);
 \coordinate (v4) at (240:1.65);
 \coordinate (v5) at (300:1.65);
 \fill[black!4] (v0)--(v1)--(v2)--(v3)--(v4)--(v5)--cycle;
 \fill[blue!12] (0,0)--(.825*1.5,.825*.8660254)--(v0)
     --(.825*1.5,-.825*.8660254)--cycle;
 \draw[thick] (v0)--(v1)--(v2)--(v3)--(v4)--(v5)--cycle;
 \draw[dashed,gray] (30:2.05)--(210:2.05);
 \draw[dashed,gray] (150:2.05)--(330:2.05);
 \node[above right] at (30:2.05) {$H_s$};
 \node[above left] at (150:2.05) {$H_t$};
 \draw[blue!55!black,thin] (0,0)--(v0);
 \fill (0,0) circle (1.4pt);
 \fill ($(v0)!.5!(v1)$) circle (1.4pt);
 \node at (.94,-.23) {$\mathcal K$};
 \foreach \i in {0,...,5} {\fill (v\i) circle (1.6pt);}
 \node[right] at (v0) {$e$};
 \node[above right] at (v1) {$s$};
 \node[above left] at (v2) {$st$};
 \node[left] at (v3) {$sts=tst$};
 \node[below left] at (v4) {$ts$};
 \node[below right] at (v5) {$t$};
\end{tikzpicture}
\caption{For $W$ of type $A_2$, $\Sigma=P_{\{s,t\}}$ is a single
hexagon with its faces. The dashed lines are the mirrors of the left
actions of $s$ and $t$; the shaded region is the Davis chamber
$\mathcal K$. The segment from $e$ to the centre divides the chamber
into two simplices of the barycentric subdivision. Boundary edges
alternate right multiplication by $s$ and $t$.}
\label{fig:davis-A2}
\end{figure}

The setwise stabilizer of $wP_J$ is $wW_Jw^{-1}$. Since $S$ is finite,
$\Sigma$ is locally finite and the action is proper and cocompact.
Moussong's theorem equips $\Sigma$ with a $W$-invariant piecewise
Euclidean CAT(0) metric \cite[Theorem~12.3.3]{Davis}: geodesic
triangles are no thicker than their Euclidean comparison triangles.
Thus any two points are joined by a unique geodesic, depending
continuously on the endpoints; the geodesics are $W$-equivariant.
In \cref{sec:geometry} these properties will give an explicit
contraction of the diagonal quotient $(\Sigma\times\Sigma)/W$.

\subsection{Based rings and the specialized Hecke algebra}\label{subsec:hecke-prelim}
A rigid based ring is a free abelian group on a distinguished basis $\mathcal B$ containing $1$, with finite-support nonnegative integral structure constants and a basis-preserving involutive anti-automorphism $b\mapsto b^{\vee}$. Its coefficient pairing satisfies
\begin{equation}\label{eq:frob}
\langle ab,c\rangle=\langle b,a^{\vee}c\rangle
=\langle a,cb^{\vee}\rangle,\qquad
\langle ab,1\rangle=\delta_{a,b^{\vee}}.
\end{equation}
Generation as a based ring allows passage to basis summands of products.

For a Coxeter system $(W,S)$, with arbitrary exponents as in
\cref{subsec:coxeter-prelim}, put
\[
 \HH(W)=\bigoplus_{w\in W}\ZZ h_w.
\]
Start with the generic Hecke algebra over $\ZZ[v,v^{-1}]$ in the normalization
$H_s^2=1+(v^{-1}-v)H_s$. Specialize by the ring homomorphism
$\ZZ[v,v^{-1}]\to\KFib$ sending $v$ to $\ph^{-1}$.
Since $\ph-\ph^{-1}=1$, the specialized generator rule is
\eqref{eq:hecke-rule}. Its structure constants are nonnegative integers,
so the $\ZZ$-span of the specialized standard basis is the subring
$\HH(W)$ displayed above. For arbitrary Coxeter systems, associativity follows from \cite[Proposition~19.1.1]{Davis}, with coefficient ring $\ZZ$ and both parameters equal to one. The reduced-word theorem gives the presentation by the quadratic relations and finite alternating braid relations. For the normalization above, see also \cite[Section~2.2]{JW}.

\begin{lemma}\label{lem:hecke}
The ring $\HH(W)$ is a rigid based ring with generator rule \eqref{eq:hecke-rule} and involution $h_w^{\vee}=h_{w^{-1}}$. It has a positive dimension function $h_w\mapsto\ph^{\ell(w)}$.
\end{lemma}
\begin{proof}
The standard Hecke recursion gives finite support and nonnegative integral coefficients. Left multiplication by $h_s$ has symmetric two-point blocks $\left(\begin{smallmatrix}0&1\\1&1\end{smallmatrix}\right)$. Hence $L_{h_w}^{\mathsf t}=L_{h_{w^{-1}}}$ and
\[
\langle h_uh_v,h_e\rangle=\langle h_v,h_{u^{-1}}\rangle.
\]
The analogous right-multiplication statement gives \eqref{eq:frob}. Reversing words yields the anti-involution. Finally, the dimension function respects the descending rule because $\ph^{k+1}=\ph^{k-1}+\ph^k$.
\end{proof}

\section{Coxeter rigidity and categorifiability restrictions}\label{sec:rigidity}
We first recover the Coxeter system from the Fibonacci relation,
without assuming a tensor categorification. Multiplication by each
generator becomes a matching with a descent term. Associativity makes the left and right
matchings commute. Their regular action then recovers a group, and the
descent indicators force Coxeter exchange. All of these steps use the rigid based-ring axioms.

\begin{theorem}[Based-ring rigidity]\label{thm:ring-rigidity}
Let $R$ be a rigid based ring generated as a based ring by distinct basis elements $x_s$, $s\in S$, where $S$ is finite and $x_s^2=1+x_s$. Then $R\cong\HH(W)$ for a Coxeter system $(W,S)$, by an isomorphism taking $x_s$ to $h_s$. Its matrix is determined by the marked two-generator subrings. At the level of based rings, the off-diagonal entries may take arbitrary values in $\{2,3,4,\ldots\}\cup\{\infty\}$.
\end{theorem}
\begin{proof}
A generator is self-dual by the unit-coefficient axiom. Its left-multiplication matrix $N$ is symmetric, nonnegative and integral, with $N^2=I+N$. At a basis vector $b$, putting $r=N_{bb}$ gives
\[
r^2+\sum_{c\ne b}N_{cb}^2=1+r.
\]
Thus $r\in\{0,1\}$ and exactly one off-diagonal coefficient is one. It defines a fixed-point-free involution $p_s$ and a function $\eps_s:\mathcal B\to\{0,1\}$ with
\[
L_s(b)=p_s(b)+\eps_s(b)b,\qquad
\eps_s(p_s(b))=1-\eps_s(b).
\]
Similarly $R_t(b)=q_t(b)+\eta_t(b)b$. Associativity gives
\begin{align*}
p_sq_t(b)-q_tp_s(b)
={}&\bigl(\eta_t(p_s(b))-\eta_t(b)\bigr)p_s(b)\\
&+\bigl(\eps_s(b)-\eps_s(q_t(b))\bigr)q_t(b).
\end{align*}
Neither basis element on the left can be $p_s(b)$ or $q_t(b)$, since that would give a fixed point of one matching. Support comparison yields $p_sq_t=q_tp_s$, and, when $p_s(b)\ne q_t(b)$,
\[
\eps_s(q_t(b))=\eps_s(b).
\]

Let $W\le\operatorname{Sym}(\mathcal B)$ be generated by the $p_s$. Based generation makes both the left and right matching actions transitive. A permutation in $W$ fixing $1$ commutes with all right matchings and hence fixes their entire orbit, so it is the identity. Identify $\mathcal B$ with $W$ using this simply transitive action. Writing $s\in W$ for the group element represented by $p_s$, we have
\[
 p_s(w)=sw,\qquad q_s(w)=ws\qquad(w\in W).
\]
Indeed, $q_s$ commutes with the left action and $q_s(e)=p_s(e)$.
The elements $s$, indexed by $S$, are distinct nonidentity involutions. The preceding direction identity becomes
\begin{equation}\label{eq:abstract-direction}
\eps_s(wt)=\eps_s(w)\quad\text{unless }sw=wt.
\end{equation}

Temporarily let $\ell$ be word length in this group. We have $\eps_s(e)=0$. If $w=t_1\cdots t_k$ is a shortest expression and $\eps_s(w)=1$, the first change of this function along its prefixes must, by \eqref{eq:abstract-direction}, occur at an index $j$ with
\[
s(t_1\cdots t_{j-1})=(t_1\cdots t_{j-1})t_j.
\]
Consequently $sw=t_1\cdots\widehat{t_j}\cdots t_k$, and the length inequality forces $\ell(sw)=k-1$. If $\eps_s(w)=0$, then $\eps_s(sw)=1$; applying the same argument to a shortest word for $sw$ gives $\ell(sw)=\ell(w)+1$. Thus the direction function is exactly the descent indicator, and the exchange condition holds for every shortest word.

Hence $(W,S)$ is a Coxeter system, with $m_{st}=\ord(st)$. The multiplication is now \eqref{eq:hecke-rule}. The two-generator subring recovers its simply transitive matching group and hence the order of $st$, proving the final assertion. Every sum in this argument is finite even when $\mathcal B$ is infinite.
\end{proof}

\begin{remark}[Torsion-freeness of the Fibonacci fusion ring]\label{rem:fib-torsion-free}
The first matrix calculation in the proof also proves that
$K_0(\Fib)$ is \emph{torsion-free}: every nonzero connected cofinite
based module is isomorphic to the regular based module, in the sense
of \cite[Definition~1.18]{ADC}. Indeed, the action
matrix $N$ of $[\tau]$ on any such module is symmetric, nonnegative
and integral, and satisfies $N^2=I+N$. The same diagonal calculation
pairs each basis vector with exactly one other; the off-diagonal
identity makes their diagonal entries sum to one. Thus $N$ is a
direct sum of blocks $\left(\begin{smallmatrix}0&1\\1&1\end{smallmatrix}\right)$.
Connectedness leaves a single block, the regular Fibonacci module.
This directly proves the torsion-freeness of the Fibonacci fusion ring
established by Arano and De Commer \cite[Proposition~1.22]{ADC}.
The block decomposition also explains the one-colour regular blocks
used below.
\end{remark}

\begin{corollary}[Descent squares]\label{lem:descent-square}
For $s,r\in S$ and $w\in W$,
\[
 sw\ne wr\quad\Longrightarrow\quad\eps_s(wr)=\eps_s(w).
\]
The analogous right-descent assertion holds. The joint orbits of $w\mapsto sw$ and $w\mapsto wr$ are consequently either two-point orbits with the same orientation, or four-point squares with parallel edges oriented consistently by length.
\end{corollary}
\begin{proof}
Apply the matching calculation in the proof of
\cref{thm:ring-rigidity} to the rigid based ring $\HH(W)$ of
\cref{lem:hecke}. Its left and right matchings are $w\mapsto sw$ and
$w\mapsto wr$. Equation~\eqref{eq:abstract-direction} gives the
displayed identity, and inversion gives its right-handed counterpart.
The matchings are commuting fixed-point-free involutions, so their
joint orbits have two or four points. On a two-point orbit they are
the same map; on a four-point orbit the descent identities give the
consistent orientations of parallel edges.
\end{proof}

The finite-composite rigidity theorem of
\cite[Theorem~\PointedRigidityTheorem]{FibPointed} has a different proof,
using Frobenius--Perron height, two-colour polygons, and local path
elimination. The proof above applies in arbitrary based rank and does not
invoke that geometric argument.

\begin{corollary}[Categorical restriction on Coxeter exponents]
\label{cor:exponent-restriction}
Let $\cD$ be a split locally finite semisimple rigid $\CC$-linear tensor
category with simple unit, generated under tensor products, finite direct
sums, and summands by finitely many pairwise nonisomorphic simple objects
$Y_s$ satisfying $Y_s^{\otimes2}\cong\one\oplus Y_s$.
Then \cref{thm:ring-rigidity} identifies $K_0(\cD)$ with $\HH(W)$ for a
marked Coxeter system $(W,S)$, and every off-diagonal exponent belongs to
$\{2,3,\infty\}$. In particular this applies when the $Y_s$ come with
specified fully faithful tensor embeddings of the positive $\Fib$.
\end{corollary}
\begin{proof}
The Grothendieck ring is a rigid based ring: semisimplicity gives
nonnegative integral coefficients of finite support, and tensor duality
gives the coefficient pairing \eqref{eq:frob}. Apply
\cref{thm:ring-rigidity}, and label the simple objects by $Y_w$, $w\in W$.
If $m_{st}=n<\infty$, the standard parabolic subgroup
$W_J$, $J=\{s,t\}$, is the finite dihedral group $I_2(n)$.

Let $\cD_J$ be the full replete additive subcategory whose simple objects
are $Y_w$ with $w\in W_J$. Every such simple is a reduced tensor word in
$Y_s,Y_t$. The Hecke generator rule keeps products with either generator
inside this list, so iteration shows closure under all tensor products in
$\cD_J$. Semisimplicity gives closure under summands. The based-ring dual
is $h_w^{\vee}=h_{w^{-1}}$, so $\cD_J$ is also closed under both tensor
duals. It contains the simple unit, has exactly $2n$ simple classes, and
inherits splitness and finite-dimensional Hom spaces from $\cD$.
It is therefore a fusion category, even if $\cD$ has infinite rank.
Its based ring is
\[
 K_0(\cD_J)\cong\HH(I_2(n))\cong R_n.
\]
The formal-codegree obstruction of
\cite[Theorem~\PointedArithmeticTheorem]{FibPointed} gives $n=2$ or $3$.
Its inputs are the formal codegrees, dimension fields and integrality
constraints of $R_n$.
\end{proof}

The exclusion depends on the fusion rules of the two generators.
The construction and uniqueness that follow fix the positive tensor
branch on every generator and establish sufficiency and marked
uniqueness for the exponents $2,3,\infty$.

\section{The six-object category of type \texorpdfstring{$A_2$}{A2}}\label{sec:A2}\label{sec:actions}
The construction in this section provides the framework for every
Coxeter system satisfying \eqref{eq:exponents}. We present the $A_2$
case first so that the actions, braid map, and reconstruction can all
be seen on six simple labels. In the general case, the additional step
is to transport the local braid maps over parabolic cosets and verify
their compatibility with every right action; this is carried out in
\cref{sec:braids}.

Concretely, let
\[
 W=W(A_2)=\{e,s,t,st,ts,o\},\qquad o=sts=tst.
\]
The aim is to turn two Fibonacci actions into a tensor category with
these six simple labels. Left actions supply the tensor generators;
right actions specify which natural transformations are admissible.
Throughout this section $S=\{s,t\}$ and $W=W(A_2)$.

\subsection{The auxiliary linear category}
Put
\begin{equation}\label{eq:M}
\cM=\Vec_W^{\mathrm{fin}}.
\end{equation}
Its objects are finitely supported families of finite-dimensional complex vector spaces; $\eobj{w}$ is the one-dimensional object in component $w$. We use $\cM$ as an auxiliary linear category on which to construct tensor actions.

Define additive linear endofunctors by
\begin{equation}\label{eq:LR}
L_s(\eobj{w})=\eobj{sw}\oplus \eobj{w}^{\oplus\eps_s(w)},\qquad
R_r(\eobj{w})=\eobj{wr}\oplus \eobj{w}^{\oplus\eta_r(w)}.
\end{equation}
For example, the two-point orbits of $L_s$, ordered by increasing
length, are $(e,s)$, $(t,st)$ and $(ts,o)$; those of $L_t$ are
$(e,t)$, $(s,ts)$ and $(st,o)$. Thus
\[
 L_s(\eobj{t})=\eobj{st},\qquad
 L_s(\eobj{st})=\eobj{t}\oplus \eobj{st},\qquad
 L_s(\eobj{o})=\eobj{ts}\oplus \eobj{o}.
\]
The right orbits are obtained by inversion. Their based modules are
regular Fibonacci modules, as explained in
\cref{rem:fib-torsion-free}. To equip the corresponding linear blocks
with full module structures, we use the fixed tensor data of $\Fib$.
On input multiplicity spaces the functors act by identity maps. On every two-point left orbit, identify the shorter vertex with $\one$ and the longer with $\tau$, and transport the complete left regular action of $\Fib$. Do the same with the complete right regular action on each right orbit. These constructions specify actions of all objects and all morphisms of the seed, including its quadratic splitting and duality. On a descent, the neighbouring vertex precedes the original vertex in the ordered direct sum.

\begin{proposition}[Elementary bimodules]\label{prop:bimodule}
For each $s,r\in S$, these left and right actions admit the structure of a $\Fib$--$\Fib$ bimodule on $\cM$.
\end{proposition}
\begin{proof}
By \cref{lem:descent-square}, a joint orbit is either two points or a consistently oriented square. On two points use the regular bimodule of $\Fib$; the middle associativity is its associator. On a square identify the underlying linear category with $\Fib\boxtimes\Fib$, with the left action on the first factor and the right action on the second. The independent-factor interchange makes it a bimodule. These identifications restrict to the already specified one-colour actions, so taking the direct sum over joint orbits proves the assertion.
\end{proof}

Write the generator interchange as
\begin{equation}\label{eq:connection}
C^0_{s,r}:L_sR_r\iso R_rL_s.
\end{equation}
A path lists its intermediate vertices in execution order; thus the source executes $R_r$ and then $L_s$. Matrix columns are source paths and rows are target paths. All one-dimensional endpoint blocks of $C^0_{s,r}$ are one. A two-dimensional block occurs exactly when
\begin{equation}\label{eq:exception}
sw=wr,\qquad\ell(sw)<\ell(w),\qquad\text{the final vertex is }w.
\end{equation}
With source and target middle vertices ordered as $(sw,w)$, that block is $U$ from \eqref{eq:U}. This gives an explicit formula for the bimodule constraint.

\subsection{The simultaneous right-module-endofunctor category}
Commutation with the right multiplication functors already lets us
recover morphisms from their value at the unit. We additionally require
compatibility with the fixed right module structures. These constraints
retain the complete seed data and will restrict the local connection
parameters in \cref{sec:local}.

Let $\cE_0$ be the category of additive linear endofunctors $T$ of $\cM$ with invertible natural module-functor constraints
\[
\gamma_{T;r,Z}:TR_{r,Z}\longrightarrow R_{r,Z}T
\]
for every right colour $r$ and every $Z\in\Fib$. They satisfy the full unit and module-functor associativity axioms, and are natural in seed morphisms. These constraints are specified separately for each right colour.

For precision, if $M\triangleleft Z$ denotes one right action and
$a^r_{M,Z,V}:(M\triangleleft Z)\triangleleft V\to M\triangleleft(Z\otimes V)$
is its associator, the required equation is
\begin{align}\label{eq:module-associativity}
a^r_{T(M),Z,V}
(\gamma_{T;M,Z}\triangleleft\id_V)
\gamma_{T;M\triangleleft Z,V}
={}&\gamma_{T;M,Z\otimes V}\,T(a^r_{M,Z,V}).
\end{align}
The colour label has been suppressed in \eqref{eq:module-associativity}; unit compatibility is imposed as well. This is the usual right-module-functor structure \cite[Section 7.2]{EGNO}.

A morphism $\eta:T\to T'$ is a natural transformation satisfying
\begin{equation}\label{eq:module-naturality}
\gamma_{T';r,Z}(\eta R_{r,Z})
=(R_{r,Z}\eta)\gamma_{T;r,Z}.
\end{equation}
Composition of these algebraic endofunctors is the tensor product, with
\begin{equation}\label{eq:gamma-product}
\gamma_{TU;r,Z}=(\gamma_{T;r,Z}U)(T\gamma_{U;r,Z}).
\end{equation}
An object of $\cE_0$ thus carries a module-functor structure for each right colour, all on the same underlying endofunctor.

\begin{corollary}[Lift of the complete seed]\label{cor:left-lift}
For every $s$, the left action is a tensor functor $\Fib\to\cE_0$, with generator $L_s$ and mixed interchanges $C^0_{s,r}$ as its right-module-functor constraints. All seed splitting and duality maps are morphisms in $\cE_0$.
\end{corollary}
\begin{proof}
The two bimodule compatibility diagrams in \cref{prop:bimodule} say exactly that the left action constraints and seed morphisms respect every right action. Hence they define the asserted tensor functor into $\cE_0$.
\end{proof}

\subsection{Evaluation at the identity}
\begin{lemma}[Evaluation propagation]\label{lem:eval-injective}
Evaluation at $\eobj{e}$ is injective on the morphisms of $\cE_0$.
\end{lemma}
\begin{proof}
Choose a reduced expression $w=r_1\cdots r_k$, and put $Q_w=R_{r_k}\cdots R_{r_1}$. All steps from $\eobj{e}$ increase length, so $Q_w(\eobj{e})=\eobj{w}$. Iterated module-functor constraints give invertible maps $\Gamma_{T,w}:TQ_w\to Q_wT$. Equation \eqref{eq:module-naturality} yields
\begin{equation}\label{eq:propagation}
\eta_{\eobj{w}}=(\Gamma_{T',w})_{\eobj{e}}^{-1}
Q_w(\eta_{\eobj{e}})(\Gamma_{T,w})_{\eobj{e}}.
\end{equation}
Only invertibility of the constraints and
\eqref{eq:module-naturality} are used; the associativity condition
\eqref{eq:module-associativity} is not needed for propagation.
The component at $\eobj{e}$ determines all simple components and,
by direct sums, all components. Thus evaluation is injective; the target Hom space is finite-dimensional.

\end{proof}

We next construct the braid isomorphism. Together with evaluation
propagation, it will identify the six simple objects and establish
semisimplicity.

\subsection{The six-object braid}
Put $A=L_sL_tL_s$ and $B=L_tL_sL_t$. Direct use of
\eqref{eq:LR} gives the following decompositions.
\begin{table}[htbp]
\centering
\caption{The common decomposition of $A(\eobj{u})$ and $B(\eobj{u})$.}
\label{tab:objects-six}
\begin{tabular}{@{}cl@{}}\toprule
$u$&Common decomposition\\\midrule
$e$&$\eobj{o}$\\
$s$&$\eobj{st}\oplus \eobj{o}$\\
$t$&$\eobj{ts}\oplus \eobj{o}$\\
$st$&$\eobj{s}\oplus \eobj{st}\oplus \eobj{ts}\oplus \eobj{o}$\\
$ts$&$\eobj{t}\oplus \eobj{st}\oplus \eobj{ts}\oplus \eobj{o}$\\
$o$&$\eobj{e}\oplus \eobj{s}\oplus \eobj{t}\oplus \eobj{st}\oplus \eobj{ts}\oplus \eobj{o}^{\oplus2}$\\\bottomrule
\end{tabular}
\end{table}
Define $T:A\to B$ to be one on every one-dimensional endpoint block. On the remaining block its source paths, in order, are
\begin{equation}\label{eq:source-braid}
(o,ts,ts,o),\qquad(o,o,o,o),
\end{equation}
and its target paths are
\begin{equation}\label{eq:target-braid}
(o,st,st,o),\qquad(o,o,o,o).
\end{equation}
Set $T=U$ on this block.

\begin{lemma}[Six-object braid]\label{lem:local-braid}
The map $T$ is an invertible module natural transformation for the complete right $s$- and $t$-actions.
\end{lemma}
\begin{proof}
Invertibility is $U^2=I$. The module structures of $A$ and $B$ already satisfy all seed axioms by \cref{cor:left-lift}. Since every seed object is a sum of $\one$ and $\tau$, it remains to check, for $r=s,t$,
\begin{equation}\label{eq:canonical-naturality}
C^0_{B,r}(TR_r)=(R_rT)C^0_{A,r}.
\end{equation}
Temporarily use
\begin{equation}\label{eq:Uz}
U(z)=\begin{pmatrix}z^2&z\\z&-z^2\end{pmatrix},\qquad P(z)=z^4+z^2-1.
\end{equation}
The six inputs and two right colours give twelve matrix equations, containing $470$ entries including structural zeros. The path computation in \cref{app:existence} gives $460$ identically zero differences over $\ZZ[z]$ and exactly ten nonzero differences, each equal to $\pm P$, $\pm zP$, or $\pm z^2P$. All ten positions and both sides are printed there. At $z=\ph^{-1/2}$ all differences vanish. Thus the displayed maps satisfy all twelve module-naturality identities.
\end{proof}

\subsection{Recovering the tensor category}
Inside $\cE_0$, take the six word functors
\[
 X_e=\Id,\quad X_s=L_s,\quad X_t=L_t,\quad
 X_{st}=L_sL_t,\quad X_{ts}=L_tL_s,\quad X_o=L_sL_tL_s.
\]
The braid just constructed identifies $X_o$ with $L_tL_sL_t$.
Their values at $\eobj{e}$ are the six distinct objects $\eobj{w}$. By
\cref{lem:eval-injective},
\[
 \Hom_{\cE_0}(X_u,X_v)=0\ (u\ne v),\qquad
 \End_{\cE_0}(X_u)=\CC.
\]
For example, $L_sX_{st}\cong X_t\oplus X_{st}$ and
$L_sX_o\cong X_{ts}\oplus X_o$. More generally, if left
multiplication increases length, concatenation gives $X_{ru}$;
if it decreases length, write $X_u\cong L_rX_{ru}$ and use
$L_r^2\cong\Id\oplus L_r$. The only change of reduced expression
needed here is $sts\leftrightarrow tst$. Thus every word in $L_s,L_t$
is a direct sum of these six functors.

Define $\Czero(A_2)$ to be their full additive subcategory in
$\cE_0$. It is tensor-closed, and the Hom calculation identifies
every endomorphism algebra with a product of full matrix algebras.
Consequently all idempotents split within this subcategory and it is
semisimple, with exactly the six simple objects $X_w$. The unit
$X_e$ is simple. The complete seed lifts of \cref{cor:left-lift}
supply duality maps for $L_s,L_t$; their tensor products and finite
direct sums supply duals for every object. This constructs a fusion
category with based ring $\HH(A_2)$, without invoking the general
reconstruction theorem.

Its dimensions are $1,\ph,\ph,\ph^2,\ph^2,\ph^3$, giving
$\FPdim\Czero(A_2)=12+16\ph$. The relation $X_sX_tX_s\cong
X_tX_sX_t$ is a Coxeter braid relation between tensor words.
Since $X_{st}\not\cong X_{ts}$, this category has no braiding.
Its identification with the known category $\TTthree$ is discussed
in \cref{sec:examples}.

\section{Actions and braid maps for general Coxeter systems}\label{sec:braids}
We now realize every Coxeter system allowed by
\cref{cor:exponent-restriction}: $S$ is finite and the off-diagonal
exponents are $2,3,\infty$, while $W$ may be infinite.
The six-object construction
separates into two steps: constructing the commuting seed actions,
and making each finite braid compatible with every right colour.
The first step uses only descent squares; the second requires
parabolic transport.

\subsection{The regular actions and their module endofunctors}
Set $\cM=\Vec_W^{\mathrm{fin}}$ and define $L_s,R_r$ by
\eqref{eq:LR}, with the same ordered summands and complete seed
actions on their two-point orbits. Define $\cE_0$ by
\eqref{eq:module-associativity}--\eqref{eq:gamma-product}, now
requiring module-functor constraints for every right colour $r\in S$.

\begin{proposition}\label{prop:general-actions}
For this arbitrary Coxeter system, the left seed actions lift to
tensor functors $\Fib\to\cE_0$, with elementary mixed interchanges
$C^0_{s,r}$ given by \eqref{eq:connection}--\eqref{eq:exception}.
Evaluation at $\eobj{e}$ is injective on morphisms of $\cE_0$.
\end{proposition}
\begin{proof}
By \cref{lem:descent-square}, every joint left--right orbit has the
same two-point or square form used in \cref{prop:bimodule}.
The regular bimodule and independent-factor bimodule therefore
give all the mixed interchanges, including compatibility with splitting
and duality maps, exactly as in \cref{cor:left-lift}.
For any $w\in W$, a reduced word supplies a finite sequence of
right actions taking $\eobj{e}$ to $\eobj{w}$. The propagation identity
\eqref{eq:propagation} determines the component of a module natural
transformation at each $\eobj{w}$ from its component at $\eobj{e}$.
Every object of $\cM$ is a finite direct sum of the $\eobj{w}$, since
it has finite support and finite-dimensional components. Additivity
of the functors and naturality with respect to the direct-sum
inclusions and projections therefore determine the transformation
on every object.
\end{proof}

\subsection{Local braid maps}
For $m_{st}=3$, use the six-object map of \cref{lem:local-braid}.
For $m_{st}=2$, use the underlying category of
$\Fib\boxtimes\Fib$, indexed by $e,s,t,st$. The interchange of
its independent left factors gives an invertible map
$L_sL_t\to L_tL_s$, with every nonzero endpoint block equal to one.
It is a module natural transformation for both local right actions.
There is no braid to impose when $m_{st}=\infty$.

\subsection{Parabolic extension}
A braid on a two-generator subsystem must commute with all right colours to become a morphism in $\cE_0$. Parabolic transport supplies this extension, including the unchanged summand on a descent.

\begin{lemma}[Parabolic transport]\label{lem:transport}
Let $d$ be shortest in $W_Jd$ and let $r\in S$. Exactly one of the following occurs:
\begin{enumerate}
\item $dr$ is a shortest representative, $W_Jd\ne W_Jdr$, and right multiplication sends $ud$ to $u(dr)$ without changing $u$. The direction of every such edge is independent of $u$.
\item There is a unique $j\in J$ with $dr=jd$. Right multiplication preserves the coset and sends $ud$ to $(uj)d$. All length comparisons agree with the corresponding comparisons in $W_J$.
\end{enumerate}
\end{lemma}
\begin{proof}
In the first case distinctness follows from uniqueness of shortest representatives and $dr\ne d$; \eqref{eq:length-add} gives the assertion about direction. Otherwise $dr$ has a left descent $j\in J$, whereas $d$ does not. By \cref{lem:descent-square}, $jd=dr$. Uniqueness follows by cancellation, and \eqref{eq:length-add} gives all remaining assertions.
\end{proof}

Fix $J=\{s,t\}$ with $m_{st}\in\{2,3\}$. On
$\Vec_{W_J}^{\mathrm{fin}}$, let $T^{\mathrm{loc}}$ denote the local
braid isomorphism: the independent-factor interchange
$L_sL_t\to L_tL_s$ when $m_{st}=2$, and the six-object braid map
$L_sL_tL_s\to L_tL_sL_t$ of \cref{lem:local-braid} when $m_{st}=3$.
On each coset $W_Jd$, with $d$ shortest, identify $\eobj{ud}$ with
$\eobj{u}$. Length additivity identifies the complete left seed
actions, including the ordered descent summands. Define $T_J$ on
this coset by replacing every vertex $v$ of a local path by $vd$
and retaining, for each pair of input and output vertices, the
corresponding path-space matrix of $T^{\mathrm{loc}}$. At exponent
three the nontrivial block is at input and output $od$.

\begin{proposition}[Global module naturality]\label{prop:global-braid}
The cosetwise maps are isomorphisms in $\cE_0$:
\[
L_sL_t\cong L_tL_s\quad(m_{st}=2),\qquad
L_sL_tL_s\cong L_tL_sL_t\quad(m_{st}=3).
\]
\end{proposition}
\begin{proof}
Let $M=(V_w)_{w\in W}$ be an object of $\cM$, and let $B_{x,w}$
denote the prescribed braid matrix between the path spaces with
input vertex $w$ and output vertex $x$. On the corresponding
summand, define $T_J$ by
\[
 B_{x,w}\otimes\id_{V_w}.
\]
A morphism $f:M\to N$ acts on these summands by the identity on
the path factor and by $f_w$ on the input multiplicity factor.
These operations commute, so $T_J$ is an ordinary natural
transformation. Its inverse is defined the same way using inverse
local matrices. It remains to prove \eqref{eq:canonical-naturality}
for \emph{every} right colour $r\in S$. Fix an input $\eobj{ud}$
and apply \cref{lem:transport}.

If $d'=dr$ is shortest, the two cosets are distinct and the right action preserves the internal coordinate $v\in W_J$. Set $\delta=\mathbf1_{\ell(d')<\ell(d)}$. Then
\[
R_r(\eobj{vd})=\eobj{vd'}\oplus \eobj{vd}^{\oplus\delta}
\quad(v\in W_J).
\]
For every left colour $q\in J$, the vertices $qvd$ and $(vd)r=vd'$ belong to different cosets. Hence the exceptional condition \eqref{eq:exception} never occurs while $R_r$ moves past the left word. All elementary interchanges simply identify the same internal left path in the two cosets. The braid matrix is identical in that coordinate before and after transport. If $\delta=1$, it is identical on the unchanged summand as well. This proves the equation in the first case.

In the other case, $dr=jd$ for a unique $j\in J$, and
\[
(vd)r=(vj)d,\qquad
\ell((vd)r)-\ell(vd)=\ell(vj)-\ell(v).
\]
Thus the entire right $r$-action becomes the local right $j$-action. This identification also preserves the connection: for $q\in J$,
\[
qvd=(vd)r\ \Longleftrightarrow\ qv=vj,
\qquad
\ell(qvd)<\ell(vd)\ \Longleftrightarrow\ \ell(qv)<\ell(v).
\]
The middle-vertex order $(qvd,vd)$ becomes $(qv,v)$, so each global elementary matrix is the local $C^0_{q,j}$. The iterated matrices and the braid therefore become precisely one of the verified local equations.

The two cases exhaust all $d$ and $r$. Unit compatibility and additivity extend the generator identity to every object of the right seed. Naturality in seed morphisms follows from the complete module structures already present. The parabolic-coordinate argument applies to every $W$ satisfying \eqref{eq:exponents}.
\end{proof}

Thus every finite braid relation is realized by a module natural
isomorphism in $\cE_0$. These maps supply the word identifications used
in the reconstruction below.

\section{Algebraic reconstruction and unitarizability}\label{sec:construction}
We now recover a tensor category from the actions. The key point is
that reduced left words are separated by evaluation at $\eobj{e}$, and
every word decomposes into reduced ones. These statements determine
the simple objects and establish semisimplicity before we pass to
duals and a positive realization.

\subsection{Normal forms before semisimplicity}
For a reduced expression $w=s_1\cdots s_k$, put
\[
X_w=L_{s_1}\cdots L_{s_k},\qquad X_e=\Id.
\]
The reduced-word theorem and \cref{prop:global-braid} identify the
objects associated with any two reduced expressions of $w$ by a
composition of module braid isomorphisms. Hence $X_w$ has a
well-defined isomorphism class.

\begin{lemma}[Separated normal forms]\label{lem:normal-forms}
One has $X_w(\eobj{e})=\eobj{w}$ and
\[
\Hom_{\cE_0}(X_w,X_v)=0\ (w\ne v),\qquad
\End_{\cE_0}(X_w)=\CC.
\]
\end{lemma}
\begin{proof}
Every suffix of a reduced word is reduced. Applying the corresponding left functors from right to left at $\eobj{e}$ therefore uses only ascending steps. The value is $\eobj{w}$. Apply \cref{prop:general-actions}; the identity of $X_w$ is nonzero.
\end{proof}

\begin{proposition}[Word decomposition]\label{prop:word-decomposition}
For every $s,w$,
\begin{equation}\label{eq:main-fusion}
L_sX_w\cong\begin{cases}
X_{sw},&\ell(sw)=\ell(w)+1,\\
X_{sw}\oplus X_w,&\ell(sw)=\ell(w)-1.
\end{cases}
\end{equation}
Every finite word in the $L_s$ is a finite direct sum of the $X_w$.
\end{proposition}
\begin{proof}
In the ascending case concatenation is reduced. In the descending case $w=s(sw)$ is reduced, and the complete seed splitting gives
\[
L_sX_w\cong L_s^2X_{sw}\cong X_{sw}\oplus L_sX_{sw}
\cong X_{sw}\oplus X_w.
\]
All these are isomorphisms in $\cE_0$. Induction proves the word decomposition.
\end{proof}

Let $\cC_{\mathrm{add}}$ be the full additive subcategory spanned by these objects. It is tensor-closed by the proposition. Define $\Czero(W)$ to be its idempotent completion; for $W=W(A_2)$ this recovers the category of \cref{sec:A2}. Concretely, objects are pairs $(T,p)$ with $p^2=p\in\End(T)$, and morphisms to $(T',q)$ satisfy $f=qfp$. The tensor product is $(T,p)(T',q)=(TT',p\otimes q)$, where the latter tensor product is that of morphisms in the endofunctor category.

\begin{proposition}[Algebraic reconstruction]\label{prop:algebraic-reconstruction}
The category $\Czero(W)$ is split locally finite semisimple and rigid, with simple unit. Its simple objects are exactly $X_w$, $w\in W$, and its based ring is $\HH(W)$.
\end{proposition}
\begin{proof}
By \cref{lem:normal-forms},
\[
\End\!\left(\bigoplus_w X_w^{\oplus n_w}\right)
\cong\prod_{w:n_w\ne0}M_{n_w}(\CC).
\]
Every idempotent in each factor is conjugate to $\operatorname{diag}(I_r,0)$. Hence every retract is isomorphic to a finite sum of the same $X_w$. This proves semisimplicity and excludes extra simple classes. Hom spaces are finite-dimensional, and $\End(X_e)=\CC$.

By \cref{prop:general-actions}, the seed evaluation and coevaluation
maps are morphisms in $\cE_0$ and satisfy the snake identities.
Hence each $L_s$ is self-dual on both sides. If $w=s_1\cdots s_k$
is reduced, taking either dual reverses the order of the factors,
giving
\[
 X_w^{\vee}\cong L_{s_k}\cdots L_{s_1}\cong X_{w^{-1}},
\]
with the analogous identification for the other dual. By the first
paragraph, every object of $\Czero(W)$ is isomorphic to a finite
direct sum of the $X_w$. Since finite direct sums of dualizable
objects are dualizable, $\Czero(W)$ is rigid.
Finally \eqref{eq:main-fusion} identifies the based ring.
\end{proof}

\begin{proof}[Proof of \cref{thm:construction}]
The preceding propositions give the category, the simple objects, the duals and the multiplication. Each lifted seed functor is fully faithful because its two simples have distinct simple images with the same scalar or zero Hom spaces. These images generate by construction. The category has finite rank exactly when $W$ is finite. In that case \cref{lem:hecke} supplies the unique positive dimension character, hence the FP dimensions in \eqref{eq:dimension-main}. For infinite $W$, the same formula $\ph^{\ell(w)}$ gives the positive dimension function on individual basis elements.
\end{proof}

\subsection{Evaluation and the reference right action}\label{subsec:reference}
\begin{corollary}\label{cor:linear-evaluation}
Evaluation at $\eobj{e}$ gives an equivalence of complex-linear categories $\Czero(W)\simeq\cM$.
\end{corollary}
\begin{proof}
It bijects the simple classes and is an isomorphism on the scalar or zero Hom spaces between them. Finite sums and idempotent ranges give the result.
\end{proof}

This linear equivalence transports the constructed tensor product, with fusion rules \eqref{eq:hecke-rule}, to $\cM$. For a module endofunctor $T$, its right-module-functor constraint at $\eobj{e}$ gives
\[
(TL_r)(\eobj{e})=T(\eobj{r})\longrightarrow R_r(T(\eobj{e})).
\]
Equation~\eqref{eq:module-associativity} makes these isomorphisms
compatible with the right action associators; naturality in seed
morphisms and unit compatibility identify the entire regular right
seed action of $\Czero(W)$ with the specified $R_r$. The right braid
morphisms of the constructed category transport along the same
identifications.

\subsection{A positive realization after construction}
\begin{proposition}[Unitarizability of the constructed category]\label{prop:unitary-realization}
The algebraic category $\Czero(W)$ admits a compatible unitary realization. In finite type it is a unitarizable fusion category. In the Hilbert realization of finite-support objects, every module natural transformation between finite sums of generating words is uniformly bounded by its component at $\eobj{e}$.
\end{proposition}
\begin{proof}
Realize $\cM$ as finitely supported graded finite-dimensional Hilbert spaces, with all linear maps as morphisms. Use the standard positive seed, orthogonal direct sums and the same path formulas. The left and right functors are $*$-functors. The seed tensor constraints, elementary interchanges and canonical braid maps are unitary, since their blocks are the identity and $U$, and $U^{\dagger}U=I$. Evaluation and coevaluation are the seed conjugate maps and satisfy the conjugate equations.

For word functors $T,T'$, taking adjoints preserves the module-naturality equation because the right-module-functor constraints are unitary. Moreover, \eqref{eq:propagation} and contractivity of the right $*$-functors give
\[
\|\eta_{\eobj{w}}\|\leq\|\eta_{\eobj{e}}\|,\qquad
\sup_{w\in W}\|\eta_{\eobj{w}}\|=\|\eta_{\eobj{e}}\|.
\]
Direct-sum naturality gives the same bound at arbitrary objects. Thus the entire algebraic Hom space between these Hilbert word functors is a finite-dimensional $C^*$-morphism space under componentwise adjoint and the displayed norm. The tensor product respects this structure.

We now split orthogonal projections in this category. Let $T$ be a
finite direct sum of generating word functors, and let
$p\in\End_{\cE_0}(T)$ satisfy $p^2=p=p^\dagger$. For each
$M\in\cM$, take the range of the orthogonal projection $p_M$ and set
\[
 T_p(M)=\im(p_M)\subseteq T(M).
\]
For a morphism $f:M\to N$, naturality gives
$T(f)p_M=p_NT(f)$, so $T(f)$ maps $T_p(M)$ into $T_p(N)$.
Its restriction defines $T_p(f)$. These restrictions make $T_p$
an additive $*$-functor. The inclusions
$i_M:T_p(M)\hookrightarrow T(M)$ form a natural isometry $i$, with
$ii^\dagger=p$ and $i^\dagger i=\id_{T_p}$.

To give $T_p$ its module-functor structure, fix a right colour $r$
and an object $Z\in\Fib$, and write
$\gamma_M=(\gamma_{T;r,Z})_M$. Since $p$ is a module natural
transformation, \eqref{eq:module-naturality} says
\[
 \gamma_M p_{R_{r,Z}(M)}
 =R_{r,Z}(p_M)\gamma_M.
\]
Thus the unitary map $\gamma_M$ carries the range of
$p_{R_{r,Z}(M)}$ onto the range of $R_{r,Z}(p_M)$. The latter is
identified isometrically with $R_{r,Z}(T_p(M))$ by
$R_{r,Z}(i_M)$, since $R_{r,Z}$ is a $*$-functor. The resulting
unitary module-functor constraint is
\[
 \gamma^p_{r,Z,M}
 =R_{r,Z}(i_M)^\dagger\,\gamma_M\,i_{R_{r,Z}(M)}
 :T_p(R_{r,Z}(M))\longrightarrow R_{r,Z}(T_p(M)).
\]
Naturality in $M$ and in seed morphisms, the unit
condition and the module-functor associativity equation all restrict
from those for $T$. Hence $T_p$ is again a simultaneous right-module
endofunctor, and $i$ is a module natural transformation.

For another projection $q$ on $T'$, the morphisms between the range
objects are identified with
\[
 \Hom(T_p,T'_q)=q\,\Hom(T,T')\,p.
\]
Taking adjoints sends this space to $p\,\Hom(T',T)\,q$, so these
spaces inherit the finite-dimensional $C^*$-category structure.
The tensor product also extends: $p\otimes q$ is an orthogonal
projection on $T\otimes T'$, and in projection-completion notation
\[
 (T,p)\otimes(T',q)=(T\otimes T',p\otimes q).
\]
By the matrix-algebra argument in \cref{prop:algebraic-reconstruction},
each range object is isomorphic to a finite direct sum of the $X_w$.
The complete seed duality maps give left and right duals for these
word objects, and finite direct sums preserve duals. The projection
completion is therefore a rigid $C^*$-tensor category.

In each of the finite-dimensional $C^*$-endomorphism algebras,
every algebraic idempotent is similar to an orthogonal projection.
Hence the projection completion is complex-linearly tensor
equivalent to the already constructed algebraic idempotent
completion. Transporting the positive structure gives the asserted
unitary realization.
\end{proof}

\subsection{Algebraic recoupling}\label{subsec:recoupling}
An $F$-matrix compares the two fusion bases for a triple tensor
product. In the endofunctor realization, composition is strictly
associative, so both bases lie in the same morphism space. This gives
an extraction procedure for any desired block.

Let $a,b,c,d\in W$ be simple-object labels in this subsection.
Choose algebraic fusion decompositions with inclusions $i_{a,b}^{c,\mu}$ and dual projections $p_{a,b}^{c,\mu}$, satisfying $pi=\delta\id$ and $\sum ip=\id$. They can be obtained by repeated use of the quadratic splittings and the braid maps. For fixed output $d$, write
\[
\lambda_{e,\mu,\nu}=(i_{a,b}^{e,\mu}\otimes\id)i_{e,c}^{d,\nu},
\qquad
\rho_{f,\sigma,\theta}=(\id\otimes i_{b,c}^{f,\theta})i_{a,f}^{d,\sigma}.
\]
If $p_\rho$ is the corresponding right-tree projection, then
\begin{equation}\label{eq:F-extraction}
p_{\rho_{f,\sigma,\theta}}\lambda_{e,\mu,\nu}
=(F_d^{a,b,c})_{(f,\sigma,\theta),(e,\mu,\nu)}\id_{X_d}.
\end{equation}
Thus columns correspond to the left-associated paths $\lambda$ and
rows to the right-associated paths $\rho$. Evaluation at $\eobj{e}$
computes each coefficient on a finite path space, even for infinite
$W$. These matrices are invertible because they compare two
decompositions of the same object. The two composites in the pentagon
compare the same initial and final fusion bases for a fourfold tensor
product, and hence agree. Orthonormal decompositions
in the positive realization give $p_\rho=\rho^{\dagger}$ and unitary
$F$-matrices. Different fusion bases change their entries by gauge
transformations while leaving the tensor category unchanged.

\section{Structural properties and finite and infinite examples}\label{sec:structure}
We examine the categories just constructed before proving their
marked uniqueness. Parabolic subsystems and products give structural
properties and examples in finite and infinite type. The $A_3$
example also relates the construction to the composition problem.

\subsection{Parabolic subcategories and products}
\begin{proposition}\label{prop:parabolic}
For $J\subseteq S$, the full tensor subcategory generated by $X_s$, $s\in J$, is marked tensor equivalent to $\Czero(W_J)$. If $S=S_1\sqcup S_2$ and $m_{st}=2$ across the partition, then
\[
\Czero(W)\simeq\Czero(W_{S_1})\boxtimes\Czero(W_{S_2}),
\]
where in infinite rank the product denotes the finite-support semisimple tensor product.
\end{proposition}
\begin{proof}
The generator rule and reduced expressions identify the parabolic simples as precisely $X_w$, $w\in W_J$. Restrict their functors to the subcategory of $\cM$ supported on $W_J$, retaining right colours in $J$. The defining actions and elementary interchanges restrict to exactly those used for $\Czero(W_J)$. Restriction is a tensor functor, bijective on simples and full on Hom spaces by evaluation, proving the first assertion.

For the second, $W=W_{S_1}\times W_{S_2}$ and length is additive. Under the corresponding factorization of $\cM$, a left or right action in one component changes only that component. Cross-component interchanges are the independent-factor ones. Sending a pair of word functors to its factorwise action is therefore a tensor functor with the specified constraints. It bijects the simples and preserves their Hom spaces, so it is an equivalence.
\end{proof}

\begin{proposition}\label{prop:subcategories}
Every full semisimple tensor subcategory of $\Czero(W)$ closed under duals and summands is generated by a subset of the marked generators.
\end{proposition}
\begin{proof}
Suppose it contains $X_w$, with $w=sv$ reduced. Then
\[
X_wX_w^{\vee}=X_sX_vX_v^{\vee}X_s
\]
contains $X_s^2=\one\oplus X_s$, so it contains $X_s$. Since $X_sX_w=X_v\oplus X_w$, it contains $X_v$. Induction extracts every generator in a reduced word for $w$. Applying this to all contained simples proves the assertion.
\end{proof}

\subsection{Finite and infinite examples}\label{sec:examples}
The finiteness criterion in \cref{subsec:geometric-representation}
identifies the finite members of the constructed family with disjoint
unions of finite ADE diagrams. We illustrate these together with an
infinite family; marked uniqueness will show that they exhaust the
categorifications with the specified diagrams.

The discrete $n$-vertex diagram gives $\Fib^{\boxtimes n}$ with $2^n$
simples. The six-object category of \cref{sec:A2} has the based ring
of $\TTthree$ in \cite[Definition~2.5]{EM}: take $X_s=\rho$,
$X_t=\mu$, $X_{st}=\overline\sigma$, $X_{ts}=\sigma$ and
$X_o=f^{(2)}$. The generator rows of that table give exactly
\eqref{eq:hecke-rule}. The unitary Fibonacci subcategories generated
by $\rho$ and $\mu$ have the positive branch. Thus
\cref{prop:rank-two} gives a marked tensor equivalence with
$\Czero(A_2)$; the fusion-table identification alone is a based-ring
identification.

A cycle with at least three vertices, with exponent three on adjacent pairs and exponent two otherwise, is a simply-laced affine diagram and gives an infinite-rank example with the same local matrices. Its objects have finite support in the infinite Coxeter group, and its structure maps are defined on finite paths. For $S=\varnothing$ the construction is simply $\Vec^{\mathrm{fd}}$.

\subsection{Free-product quotients}
The positive realization gives a free-product description of the constructed categories. This will be used for the three-generator example in \cref{subsec:A3}.

\begin{proposition}[The free-product quotient]\label{prop:free-product-quotient}
For nonempty $S$, the specified seed functors extend, up to monoidal natural isomorphism on each factor, to a faithful dominant strong tensor functor
\[
 Q:\mathop{*}_{s\in S}\Fib\longrightarrow\Czero(W).
\]
It can be chosen unitary in the realization of \cref{prop:unitary-realization}. The images of the nontrivial objects of the specified factors are pairwise nonisomorphic. For $S=\varnothing$, the corresponding functor is the identity of $\Vec^{\mathrm{fd}}$.
\end{proposition}
\begin{proof}
In the positive realization the full seed functors are unitary. The universal property for free products of rigid $C^*$-tensor categories \cite[Theorem~B]{Valvekens} extends them to a unitary tensor functor. Transport it along the equivalence with the algebraic model to obtain $Q$.

We verify faithfulness explicitly. For every nonzero simple $Z$ in the source, rigidity gives a nonzero unit-channel map and semisimplicity makes $\one$ a direct summand of $Z^{\vee}\otimes Z$. If $Q(Z)=0$, additivity and strong tensoriality would then imply $Q(\one)=0$, contrary to the unit isomorphism. Hence $Q(Z)\ne0$, and the map $\CC\id_Z\to\End(Q(Z))$ is injective. Decompose two arbitrary objects into finite sums of simples. Composing a morphism with the source projections and inclusions extracts its scalar matrix entries; the same extraction after applying $Q$ proves that a morphism killed by $Q$ has every entry zero. 

If $w=s_1\cdots s_k$ is reduced, the generator rule gives $X_{s_1}\otimes\cdots\otimes X_{s_k}\cong X_w$. Every target simple therefore lies in the image of a source word, proving dominance. Finally, distinct generators give distinct standard basis labels $h_s$, proving the last assertion for nonempty $S$.
\end{proof}

\subsection{The category of type \texorpdfstring{$A_3$}{A3}}\label{subsec:A3}\label{sec:associator-counts}\label{subsec:A3-direct}
Our construction of $\Czero(A_3)$ in
\cref{sec:braids,sec:construction} is independent of direct
$F$-symbol computation. The calculation below illustrates the usual
approach through pentagon equations and fusion-space basis choices.

For the three-vertex path $s$--$t$--$u$, the group $W(A_3)=S_4$ labels the $24$ simple objects. Its length enumerator is
\[
1+3t+5t^2+6t^3+5t^4+3t^5+t^6.
\]
Substitution $t=\ph^2$ gives
\[
 \FPdim\Czero(A_3)=300+480\ph=540+240\sqrt5.
\]
The size of a full associator table is determined by these fusion rules.
Let $N_{ab}^{x}$ be the multiplicity of $X_x$ in $X_a\otimes X_b$, and set
\[
 m_{abc}^{d}=\sum_{x\in W}N_{ab}^{x}N_{xc}^{d}
            =\sum_{y\in W}N_{bc}^{y}N_{ay}^{d}.
\]
For each ordered quadruple $(a,b,c,d)$ with $m_{abc}^{d}>0$, the
associator has an $m_{abc}^{d}\times m_{abc}^{d}$ matrix block. In type
$A_3$, the numbers of blocks and scalar matrix positions are
\[
 \#\{(a,b,c,d):m_{abc}^{d}>0\}=126{,}802,\qquad
 \sum_{a,b,c,d\in W}(m_{abc}^{d})^2=1{,}293{,}175.
\]
The largest block has size $100\times100$. These counts include
quadruples involving the unit object and count positions with
multiplicity, before imposing gauge choices or coherence equations.
The script \path{A3_Associator_Counts.py} computes the fusion products
from \eqref{eq:hecke-rule} and checks the two expressions for every
$m_{abc}^{d}$ using exact integers.

For this direct calculation we use nonorthogonal algebraic fusion
bases. Put $d=\ph$ and rescale the unit-channel projection and inclusion
of each Fibonacci seed by
\[
 p'_{\one}=\sqrt d\,p_{\one},\qquad
 i'_{\one}=d^{-1/2}i_{\one},
\]
leaving the $\tau$-channel and strict-unit maps unchanged. In the
channel order $(\one,\tau)$, the associator in \eqref{eq:U} becomes
\[
 DUD^{-1}=\begin{pmatrix}d-1&1\\d-1&1-d\end{pmatrix},
 \qquad D=\operatorname{diag}(\sqrt d,1).
\]
The evaluation and coevaluation of the same abstract seed are now
$p'_{\one}$ and $d\,i'_{\one}$; the snake identities reduce to
$d(d-1)=1$. The two adjacent $\TTthree$ models are expressed in
compatible algebraic fusion bases extending this normalization on
their common Fibonacci seed. For the commuting
$\Fib\boxtimes\Fib$ model we use the tensor-product bases of the
same normalized seeds. Thus the three local models agree on their
overlaps. The fixed local boundary in the calculation below refers
to these rationalized local tables.

The first
reduction keeps only associators whose third input is a generator.
Indeed, whenever a chosen reduced word gives $X_k\otimes X_l\cong X_v$
with $l$ a generator, the pentagon expresses the associator with
third input $X_v$ through those with shorter third inputs. This
reduces the calculation to a core table of $8{,}316$ blocks and $15{,}660$ positions.
Comparing the reconstructions along $sts=tst$ and $tut=utu$ imposes
additional equations on that core.

The next reduction uses the full matrix freedom in fusion-space
bases. Along each chosen reduced word, bases can be selected
recursively so that the corresponding prefix-extension associator
is the identity, except on the fixed local boundary. Word length orders these compatible choices. Linear and monomial
elimination, dividing only by quantities proved nonzero, then leaves
ten nonzero coordinates. A residual scalar gauge removes all ten:
its integer weight matrix has a maximal minor of determinant $1$,
so this residual normalization uses integer powers and requires no
additional root choices.

In this algebraic gauge, reconstruction yields the full table over
$\ZZ[d]$, where $d^2=d+1$;
the positive embedding is $d=\ph$. All pentagons were checked in
exact arithmetic, and all associator blocks were checked to be
invertible.

The same category also gives a concrete comparison with the
three-generator quotient problem.

\begin{proposition}[Type $A_3$ and a three-generator quotient]\label{prop:A3-comparison}
Let $s=(12)$, $t=(23)$ and $u=(34)$ in $W(A_3)=S_4$, and put $a=X_s$, $b=X_t$, $c=X_u$. This gives a rank-preserving finite unitarizable quotient of $\Fib^{*3}$ in the sense of \cite[Definitions~2.3, 2.4 and 2.6]{EM}. It satisfies
\[
 ac\cong ca,\qquad aba\cong bab,\qquad bcb\cong cbc,
 \qquad ab\not\cong ba,\quad bc\not\cong cb.
\]
Consequently each of $a,b,c$ fails to commute with at least one of the
other two generators, contrary to \cite[Lemma~4.7, p.~819]{EM}.
The commuting pattern and the length-five relation below also contradict
\cite[Sublemma~4.7.4, pp.~830--831]{EM}.
\end{proposition}
\begin{proof}
The construction gives a rigid category with $24$ simple classes, \cref{prop:unitary-realization} gives unitarizability, and \cref{prop:free-product-quotient} gives the faithful dominant tensor functor. The distinct simple objects $a,b,c$ verify rank preservation, which in the cited definition means precisely distinctness of these designated generators. The displayed isomorphisms follow from reduced Coxeter words. The products $st,ts$ are distinct, as are $tu,ut$, so the corresponding simple objects are nonisomorphic.

We spell out the short relation relevant to the proof of the cited lemma:
\begin{equation}\label{eq:A3-five}
 abcab\cong abacb\cong babcb\cong bacbc\cong bcabc.
\end{equation}
The four steps use, in order, the commuting $a,c$ pair, the $a,b$ braid, the $b,c$ braid, and the commuting $a,c$ pair. Let $v_n^a,v_n^b,v_n^c$ be the length-$n$ prefixes of $(abc)^\infty,(bca)^\infty,(cab)^\infty$. With permutations composed from right to left, their one-line labels are
\begin{equation}\label{eq:A3-prefixes}
\begin{array}{c|ccc}
n&v_n^a&v_n^b&v_n^c\\\hline
1&2134&1324&1243\\
2&2314&1342&2143\\
3&2341&3142&2413\\
4&3241&3412&2431\\
5&3421&3421&4231
\end{array}
\end{equation}
Each permutation in row $n$ has $n$ inversions. Thus these are reduced words and their products are simple. Their labels are pairwise distinct for $n<5$, and the first two coincide for $n=5$; every word in \eqref{eq:A3-five} has label $3421$ and length five. Reversing words preserves simplicity by duality. This verifies the minimal length and the simplicity assumptions used in the cited proof. In its case where only $a$ and $c$ commute, \cite[Sublemma~4.7.4, pp.~830--831]{EM} asserts exclusion of the small cases $n=1,\ldots,5$; it works out $n=2$ and leaves the other cases to the reader. Relation \eqref{eq:A3-five} shows that this exclusion cannot hold.
\end{proof}

A tensor autoequivalence carrying $a\mapsto b\mapsto c\mapsto a$
would send $ac\cong ca$ to $ba\cong ab$, contradicting the distinct
simple labels $ts$ and $st$. Thus the $A_3$ quotient lies outside the
cyclic class of \cite[Definition~2.7 and Theorem~4.1]{EM}.
The published argument uses Lemma~4.7 in the proof of Lemma~4.8
(pp.~831--832) to obtain a tensor-product block decomposition, then
imposes cyclic symmetry in Remark~4.9 and Corollary~4.10
(pp.~832--833). The failure of Lemma~4.7 invalidates that general
block-decomposition step. A proof of the cyclic classification must
establish the required decomposition under the cyclic hypothesis;
the $A_3$ example addresses the broader intermediate claim.

\section{Regular module actions and the tensor-equivalence problem}\label{sec:regular}
The construction supplies a categorification for every allowed
Coxeter system. We now ask whether the same marked fusion rules admit
any other tensor structure. Let $\cD$ be a marked categorification as
in \cref{def:marked}, with simple objects $Y_w$ labelled by $h_w$.
The reference category $\Czero(W)$ of \cref{sec:construction} has
simples $X_w$, associator $a^0$, and marked seed functors
$\iota_s^0:\Fib\to\Czero(W)$.
Their common simple labels give a linear equivalence; the problem is to
make it a tensor equivalence preserving the marked Fibonacci functors.
We first describe the regular module blocks inside $\cD$ and their identifications,
and then explain what their simultaneous compatibility must accomplish.

\subsection{The tensor equivalence to be constructed}
\label{subsec:regular-target}
Choose a label-preserving $\CC$-linear equivalence
\[
E:\Czero(W)\longrightarrow\cD,\qquad E(X_w)=Y_w.
\]
To make $E$ a strong tensor functor, one needs natural isomorphisms
\begin{equation}\label{eq:tensorator-target}
J_{A,B}:E(A)\otimes E(B)\iso E(A\otimes B)
\end{equation}
and a unit isomorphism, satisfying
\begin{equation}\label{eq:tensorator-coherence}
\begin{aligned}
&J_{A,B\otimes C}(\id\otimes J_{B,C})
 a^{\cD}_{E(A),E(B),E(C)}\\
&\qquad =E(a^0_{A,B,C})J_{A\otimes B,C}(J_{A,B}\otimes\id).
\end{aligned}
\end{equation}
Our convention is $a_{A,B,C}:(A\otimes B)\otimes C\to A\otimes(B\otimes C)$, and compositions are read from right to left. A marked tensor equivalence also carries monoidal natural isomorphisms $E\iota_s^0\cong\iota_s$ for all $s$.

The local choices used below identify the left and right regular blocks of each marked seed with fixed models. Each block identification is a module equivalence, including its structure isomorphisms. These identifications are chosen separately; they are not yet prescribed components of a common family \eqref{eq:tensorator-target}. The proof first chooses them so that the mixed interchanges and normalized braid matrices agree with the reference, and then obtains the full tensorator from regular realizations.

\subsection{Regular blocks inside the candidate category}
\label{subsec:intrinsic}
Fix $s\in S$ and $x\in W$ with $y=sx$ and $\ell(y)=\ell(x)+1$. Let $\cD_{s,x}$ be the full additive subcategory on $Y_x,Y_y$. The restricted action of $\Fib$ on $\cD$ makes this a left module category. The functor
\begin{equation}\label{eq:canonical-phi}
\Phi^L_{s,x}:\Fib\longrightarrow\cD_{s,x},\qquad
V\longmapsto\iota_s(V)\otimes Y_x,
\end{equation}
is a module equivalence. Indeed, it sends the two simple objects to $Y_x$ and an object isomorphic to $Y_y$, and it is fully faithful on their scalar or zero morphism spaces. Its module constraint is induced by $a^{\cD}$ and the tensorator of $\iota_s$. The pentagon in $\cD$ gives the module-functor axiom. The analogous right equivalence is
\begin{equation}\label{eq:canonical-phi-right}
\Phi^R_{r,x}:\Fib\longrightarrow\cD_{x,r},\qquad
V\longmapsto Y_x\otimes\iota_r(V),
\quad \ell(xr)=\ell(x)+1.
\end{equation}
Here $\cD_{x,r}$ is the full additive subcategory on $Y_x,Y_{xr}$.

To express \eqref{eq:canonical-phi} in the fixed simple-object labels, choose an additive linear equivalence
\[
K_{s,x}:\Fib\longrightarrow\cD_{s,x},\qquad
K_{s,x}(\mathbf1)=Y_x,\quad K_{s,x}(\tau)=Y_y,
\]
and a linear natural isomorphism $\theta:K_{s,x}\Rightarrow\Phi^L_{s,x}$. Its component at $\mathbf1$ is fixed by the unit identifications. Its component at $\tau$ is a choice of nonzero map
\begin{equation}\label{eq:ascending-inclusion}
\theta_\tau:Y_y\iso Y_s\otimes Y_x.
\end{equation}
Transport the module structure of $\Phi^L_{s,x}$ through $\theta$. If
\[
\phi_{Z,V}:\Phi^L_{s,x}(Z\otimes V)
   \longrightarrow\iota_s(Z)\otimes\Phi^L_{s,x}(V)
\]
is its fixed module constraint, the transported constraint on $K_{s,x}$ is
\begin{equation}\label{eq:transport-kappa}
\kappa_{Z,V}
=(\id_{\iota_s(Z)}\otimes\theta_V^{-1})\,
 \phi_{Z,V}\,\theta_{Z\otimes V}.
\end{equation}
This specifies a complete module equivalence, including compatibility with all seed morphisms.

The map \eqref{eq:canonical-phi} and its module structure remain fixed. The choice in \eqref{eq:ascending-inclusion} does not: replacing $\theta_\tau$ by $p\theta_\tau$, for $p\in\Cx$, produces another identification \eqref{eq:transport-kappa}. The corresponding changes on the upward, downward and stationary channels are $p,p^{-1},1$, respectively. Thus the freedom is already present in $\cD$: it consists of choosing the maps that identify its regular module blocks with labelled copies of the fixed regular seed. \Cref{lem:edge-gauge} will show that these are exactly the possible changes.

\subsection{A common regular realization}
\label{subsec:regular-realization}
We now place all the labelled blocks in a common linear category
\[
\cM=\Vec^{\mathrm{fin}}_W.
\]
Its simple object in degree $w$ is denoted $\mathbf e_w$. Choose a label-preserving linear equivalence $F:\cD\to\cM$ and a quasi-inverse. The transported seed actions are
\[
\widetilde L_{s,Z}=F(\iota_s(Z)\otimes-)F^{-1},\qquad
\widetilde R_{r,Z}=F(-\otimes\iota_r(Z))F^{-1}.
\]
The block identifications above compare each of them with a fixed complete regular action $L_{s,Z}$ or $R_{r,Z}$. On the generator their ordered decompositions are \eqref{eq:LR}.
Each two-point block carries the full regular seed action, including
its associativity constraints and seed morphisms.

In these coordinates, the mixed interchange is
\begin{equation}\label{eq:mixed-C}
C_{s,r}:L_sR_r\iso R_rL_s.
\end{equation}
It is the coordinate form of
\[
a_{Y_s,M,Y_r}^{\cD,-1}:
Y_s\otimes(M\otimes Y_r)\longrightarrow(Y_s\otimes M)\otimes Y_r.
\]
Its compatibility with the complete left and right actions follows from the bimodule axioms inherited from $\cD$. As a right-module-functor constraint for $L_s$, it is the component $\gamma_{L_s;r,\tau}$; the term mixed interchange emphasizes its role between the two actions.

Recall from \cref{sec:braids} that $\cE_0$ consists of additive
linear endofunctors of $\cM$ equipped, separately for each right
colour, with a complete right-module-functor structure. Its objects
satisfy \eqref{eq:module-associativity}, its morphisms satisfy
\eqref{eq:module-naturality}, and composition is its tensor product,
as in \eqref{eq:gamma-product}. No compatibility between distinct
right colours is imposed in this definition.

Left tensoring gives a strong monoidal functor
\begin{equation}\label{eq:regular-functor}
\operatorname{Reg}_{\cD}:\cD\longrightarrow\cE_0,\qquad
X\longmapsto F(X\otimes-)F^{-1}.
\end{equation}
The chosen right block identifications transport its module constraints to the fixed right actions. Its tensorator comes from $a^{\cD}$ and the equivalence data for $F$. The chosen left block identifications present each marked generator as $L_{s,Z}$ with its transported right-module-functor constraints. The reference category has the same kind of realization, denoted $\operatorname{Reg}_0$.

\begin{lemma}[Full faithfulness of the regular realization]
\label{lem:evaluation}
The strong monoidal functor $\operatorname{Reg}_{\cD}$ is fully
faithful.  A module natural transformation between two functors in its
image is determined by its component at $\eobj{e}$, with the propagation
formula \eqref{eq:propagation}.
\end{lemma}

\begin{proof}
The evaluation argument of \cref{prop:general-actions} shows that a module natural
transformation in $\cE_0$ is determined by its value at $\eobj{e}$:
for a reduced right word representing $w$, the corresponding invertible
right-module constraints give precisely \eqref{eq:propagation}.  This proves faithfulness.

Let $X,Y\in\cD$.  Under the unit identifications,
$\operatorname{Reg}_{\cD}(X)(\eobj{e})\cong F(X)$ and similarly for $Y$.
Thus any component at $\eobj{e}$ of a transformation between their regular
functors corresponds to a unique morphism $f:X\to Y$.
The transformations $f\otimes\id_{(-)}$, transported by $F$,
are right-module natural by naturality of the associator; they have the
prescribed value at $\eobj{e}$.  Propagation shows there are no others.
Hence the functor is full as well as faithful.
\end{proof}

The same full-faithfulness argument applies in the larger category
where only natural commutation isomorphisms with the right multiplication
functors are retained, with morphisms satisfying
\eqref{eq:module-naturality}. More generally, for any split locally
finite semisimple categorification, commutation isomorphisms with all simple right
multipliers determine a transformation between left tensoring functors
from its value $f$ at the unit: its component at $Y$ is
$f\otimes\id_Y$, under the unit identifications. This is the universal
regular realization described in the Introduction.

\subsection{Freedom in the one-colour block identifications}
\label{subsec:one-colour}
Fix a left seed action $Z\mapsto L_Z$ on $\cM$, with tensor structure $\mu_{Z,V}:L_ZL_V\to L_{Z\otimes V}$. Two label-preserving identifications with this action differ, after a linear natural identification, by a module autoequivalence with underlying functor $\Id_{\cM}$. Its structure is a family
\[
\gamma_{Z,M}:L_Z(M)\iso L_Z(M)
\]
natural in both variables, equal to the identity on the unit, and satisfying
\begin{equation}\label{eq:left-comparison-coherence}
\gamma_{Z\otimes V}\mu_{Z,V}
=\mu_{Z,V}(\gamma_Z*\gamma_V),\qquad
(\gamma_Z*\gamma_V)_M=\gamma_{Z,L_V(M)}L_Z(\gamma_{V,M}).
\end{equation}
In particular, naturality in $Z$ includes the chosen seed splitting and duality maps. The family $\gamma$ is the structure of a module functor; an isomorphism between two such module functors is a further datum.

\begin{lemma}[Rescaling a regular module block]\label{lem:edge-gauge}
On a regular two-point block $\{x,y\}$, with $x$ shorter, every module-functor structure on the identity has
\begin{equation}\label{eq:edge-gauge}
\gamma_{\tau,\mathbf e_x}=p\,\id_{\mathbf e_y},\qquad
\gamma_{\tau,\mathbf e_y}=\operatorname{diag}(p^{-1},1)_{\mathbf e_x\oplus\mathbf e_y},
\qquad p\in\Cx.
\end{equation}
Every $p$ occurs, and the resulting module autoequivalence is module naturally isomorphic to the identity. The same statements hold for right regular actions.
\end{lemma}
\begin{proof}
Write the three scalar coefficients of $\gamma_\tau$ as $p,q,t$, so that its components are $p$ and $\operatorname{diag}(q,t)$. Apply \eqref{eq:left-comparison-coherence} at $Z=V=\tau$, $M=\mathbf e_x$. On the unit and $\tau$ channels, naturality in the seed inclusions makes the left side act by $1,p$, whereas the right side acts by $pq,pt$. Hence $pq=1$ and $pt=p$, giving $q=p^{-1}$ and $t=1$.

Conversely, let $\eta$ be the natural automorphism of the identity linear functor with $\eta_{\mathbf e_x}=1$, $\eta_{\mathbf e_y}=p$. Define
\begin{equation}\label{eq:coboundary-gamma}
\gamma_{Z,M}=\eta_{L_Z(M)}L_Z(\eta_M^{-1}).
\end{equation}
Transport of the trivial module structure through $\eta$ proves all coherence and naturality axioms. The coefficients are \eqref{eq:edge-gauge}, and
$L_Z(\eta_M)\gamma_{Z,M}=\eta_{L_Z(M)}$ makes $\eta$ a module natural isomorphism to the identity module functor. The right-module argument is identical.
\end{proof}

For each colour we may make these choices independently on its disjoint regular blocks. We denote the upward left and right coefficients by
\begin{equation}\label{eq:edge-scalars}
A_s(x)\quad(\ell(sx)>\ell(x)),\qquad
B_r(x)\quad(\ell(xr)>\ell(x)).
\end{equation}
We call these changes \emph{edge rescalings}. By
\cref{lem:edge-gauge}, each upward scalar specifies a complete
module-functor structure, with its inverse on the downward channel
and $1$ on the stationary channel.
The same coefficient is then used whenever that one-colour block occurs in a longer word. Different colours, and left and right occurrences, retain separate block identifications. Even when two occurrences have the same endpoints, their identifications need not agree. A single global natural automorphism of $\Id_{\cM}$ would instead give vertex scalars and only endpoint ratios on all incident edges. That smaller freedom has trivial product around every closed edge path.

\subsection{Mixed associativity in fusion coordinates}
\label{subsec:path-change}
The effect of \eqref{eq:edge-scalars} can be computed directly in $\cD$. Consider a four-point orbit with shortest vertex $x$ and vertices $x,sx,xr,sxr$. Write
\[
\ell_{s,x}:Y_{sx}\iso Y_s\otimes Y_x,\qquad
r_{r,x}:Y_{xr}\iso Y_x\otimes Y_r
\]
for the ascending identifications coming from \eqref{eq:ascending-inclusion} and its right analogue. The two path inclusions are
\begin{align*}
I_\alpha&=(\id_{Y_s}\otimes r_{r,x})\ell_{s,xr}:
Y_{sxr}\longrightarrow Y_s\otimes(Y_x\otimes Y_r),\\
I_\beta&=(\ell_{s,x}\otimes\id_{Y_r})r_{r,sx}:
Y_{sxr}\longrightarrow(Y_s\otimes Y_x)\otimes Y_r.
\end{align*}
\[
\begin{tikzpicture}[>=Stealth,every node/.style={font=\small}]
 \node (x) at (0,0) {$x$};
 \node (sx) at (0,1.8) {$sx$};
 \node (xr) at (4,0) {$xr$};
 \node (sxr) at (4,1.8) {$sxr$};
 \draw[->] (x)--node[left] {$L_s:\ A_s(x)$}(sx);
 \draw[->] (xr)--node[right] {$L_s:\ A_s(xr)$}(sxr);
 \draw[->] (x)--node[below] {$R_r:\ B_r(x)$}(xr);
 \draw[->] (sx)--node[above] {$R_r:\ B_r(sx)$}(sxr);
\end{tikzpicture}
\]
There is a unique $h\in\Cx$ such that
\begin{equation}\label{eq:intrinsic-square}
a^{\cD,-1}_{Y_s,Y_x,Y_r}I_\alpha=I_\beta h.
\end{equation}
Thus $h$ is a matrix coefficient of the fixed associator in the chosen fusion coordinates. Under $\ell'_{s,x}=\ell_{s,x}A_s(x)$ and $r'_{r,x}=r_{r,x}B_r(x)$, it becomes
\begin{equation}\label{eq:square-change-intro}
h'=h\,\frac{B_r(x)A_s(xr)}{A_s(x)B_r(sx)}.
\end{equation}
Both the tensor category and its associator are unchanged.

More generally, for a path $\alpha$ through left and right actions, let $g(\alpha)$ be the product of its edge coefficients. A downward traversal contributes the inverse upward coefficient, and a stationary descent channel contributes $1$. For its inclusion and projection,
\[
i'_\alpha=i_\alpha g(\alpha),\qquad
p'_\alpha=g(\alpha)^{-1}p_\alpha.
\]
Consequently a matrix $M$ of a natural transformation changes by
\begin{equation}\label{eq:path-gauge}
M'_{\beta,\alpha}=g(\beta)^{-1}M_{\beta,\alpha}g(\alpha).
\end{equation}
This formula applies to mixed interchanges and braid maps, with
edge coefficients restricted to the module equivalences of
\cref{lem:edge-gauge}.

\subsection{The marked seed has no residual tensor symmetry}

The markings specify full Fibonacci tensor functors. The following
fact will identify these functors when constructing the local bimodule equivalences.

\begin{lemma}[Tensor autoequivalences of the Fibonacci seed]
\label{lem:fib-aut}
Every $\CC$-linear tensor autoequivalence of $\Fib$
which fixes the isomorphism class of $\tau$ is monoidally naturally
isomorphic to the identity tensor functor.
\end{lemma}

\begin{proof}
After a linear natural identification, assume that the functor is the
identity on objects and morphisms, with strict unit.  Its tensorator
on $\tau\otimes\tau\cong\one\oplus\tau$ is diagonal,
$\operatorname{diag}(\lambda,\mu)$, where
$\lambda,\mu\in\Cx$.
The triple-$\tau$ monoidal coherence condition, on the output-$\tau$
space with intermediate channels $(\one,\tau)$, implies
\[
 \operatorname{diag}(\lambda,\mu^2)\,U
   =U\,\operatorname{diag}(\lambda,\mu^2).
\]
The off-diagonal entries of $U$ are nonzero, hence
$\lambda=\mu^2$.  The natural automorphism which is $1$ on
$\one$ and $\mu$ on $\tau$ trivializes the tensorator.
\end{proof}

\subsection{Recovering a tensor equivalence from matching actions}
\label{subsec:formal-reconstruction}
The following observation separates the categorical reconstruction from the scalar normalization problem.

\begin{proposition}[A criterion for marked tensor equivalence]\label{prop:comparison-criterion}
Suppose that changes of the complete left and right one-colour identifications make all candidate mixed interchanges equal to the reference interchanges. Then $\cD$ and $\Czero(W)$ are marked tensor equivalent.
\end{proposition}
\begin{proof}
Let $\beta_{r,Z}:R_{r,Z}\to R_{r,Z}$ be the right coordinate changes. Because each $\beta_r$ is the constraint family of a right-module autoequivalence with underlying functor $\Id_{\cM}$, the identity linear functor with constraints $\beta_{r,Z}^{-1}$ is an invertible object $Q$ of $\cE_0$. There is no compatibility condition between different $r$ in the definition of $\cE_0$. Conjugation $T\mapsto QTQ^{-1}$ is a monoidal autoequivalence and changes the module-functor constraint of $(T,\gamma)$ to
\begin{equation}\label{eq:conjugation-formula}
\gamma'_{T;r,Z}
=(\beta_{r,Z}T)^{-1}\gamma_{T;r,Z}(T\beta_{r,Z}).
\end{equation}
This follows directly from \eqref{eq:gamma-product} and realizes the right coordinate change.

After this conjugation, the left block identifications give monoidal natural isomorphisms between the candidate and reference marked actions as tensor functors $\Fib\to\cE_0$. Their compatibility with the right-module structures is exactly the equality of the transported mixed interchanges. Equality for the pair of seed generators extends to arbitrary seed
objects in both variables by additivity and naturality: each object of $\Fib$ is a sum of $\one$ and $\tau$,
and the block identifications already preserve all seed morphisms. Thus the two
regular realizations have the same marked generating tensor functors
up to monoidal natural isomorphism.

Both realizations are fully faithful by \cref{lem:evaluation}. Their replete essential images are the full tensor subcategories generated by these actions under finite sums and summands. Indeed, full faithfulness lifts an idempotent on an image object to an idempotent in the source, where it splits. The two images therefore coincide. Composing the equivalences onto this common image gives a tensor equivalence. It preserves the labels, because a reduced word for $w$ represents the simple object with label $w$, and it preserves the markings through these monoidal natural isomorphisms on the seed functors.
\end{proof}

It remains to choose the edge rescalings in \eqref{eq:edge-scalars} so that all mixed interchanges are normalized simultaneously. The same family will also normalize the coordinate parameters of the intrinsic braid morphisms. These extra parameters supply the pure faces needed for the cellular argument.

\section{Local bimodule structures}\label{sec:local}
The regular realization reduces the comparison of candidates to their
mixed interchanges. We first identify the bimodule category on each
joint orbit, and then express its interchange in the fixed one-colour
coordinates of \cref{sec:regular}. There are two local models: the
independent left and right regular actions on four simples, and the
regular Fibonacci bimodule on two simples. Their coordinate freedoms
are described by one nonzero scalar in each case.

After the complete actions have been fixed, write
\[
 C_{s,r}:L_sR_r\longrightarrow R_rL_s
\]
for a candidate interchange. We describe all possibilities on the joint orbits of \cref{lem:descent-square}. Matrix columns are source paths and rows are target paths. Paths list vertices in execution order, so the source operation list here is $(R_r,L_s)$ and the target list is $(L_s,R_r)$.

\begin{lemma}[Local scalar classification]\label{lem:local-phases}
Every candidate mixed interchange on a joint orbit belongs to one of the following complete families.
\begin{enumerate}
\item On a four-point orbit $x,sx,xr,sxr$, with $x$ shortest, there is one scalar $h\in\Cx$. The channels on which both operations move to the neighbouring vertex have coefficients $h,h^{-1},h^{-1},h$ at those four inputs, respectively. All other scalar blocks are $1$.
\item On a two-point orbit $x,y=sx=xr$, with $x$ shorter, there is one scalar $k\in\Cx$. The scalar blocks with input/output pairs $(x,x),(x,y),(y,x)$ are $k^2,k,k$, respectively. The $(y,y)$ block, with middle vertices ordered as $(x,y)$, is
\begin{equation}\label{eq:Uk}
 U_k=\begin{pmatrix}ak^{-2}&bk^{-1}\\bk^{-1}&-a\end{pmatrix}
 =\operatorname{diag}(k^{-1},1)U\operatorname{diag}(k^{-1},1).
\end{equation}
\end{enumerate}
Every value of the displayed parameter gives a genuine bimodule structure on the indicated joint orbit.
\end{lemma}
\begin{proof}
On a four-point orbit, let $\cM_{s,r,x}$ be the full additive
subcategory on $Y_x,Y_{sx},Y_{xr},Y_{sxr}$. The functor
\[
 \Fib\boxtimes\Fib\longrightarrow\cM_{s,r,x},
 \qquad (A,B)\longmapsto\iota_s(A)\otimes Y_x\otimes\iota_r(B),
\]
intertwines the left action on the first factor and the right action
on the second. Its module-functor constraints come from the candidate
associator and the tensorators of the seed embeddings; their bimodule
compatibility follows from the pentagon. The four simple pairs map to
the four distinct orbit simples, so the underlying linear functor is
fully faithful and essentially surjective. This gives a bimodule
equivalence with the independent-factor model.

On a two-point orbit, first identify the right action with the regular
right $\Fib$-module. Evaluation at $\one$ gives a monoidal equivalence
\[
 \operatorname{End}_{\Fib}(\Fib)\simeq\Fib,
 \qquad T\longmapsto T(\one),
\]
whose inverse sends $A$ to left tensoring by $A$. Indeed, the
right-module constraint identifies $T(Z)$ with $T(\one)\otimes Z$;
its compatibility with the module associator makes this a module
natural identification. Composition then corresponds to tensor product.
The commuting left action is a tensor functor into this endofunctor
category. Under evaluation it becomes a tensor endofunctor of $\Fib$
sending $\tau$ to $\tau$. It is an equivalence on the two simple
objects, and \cref{lem:fib-aut} identifies it monoidally with the
identity. Thus the entire two-point bimodule is equivalent to the
regular Fibonacci bimodule.

These bimodule equivalences supply block identifications in which the
interchanges are the reference ones. We now express them in the
previously chosen one-colour coordinates. Identify the underlying
labelled linear equivalences with the identity. By
\cref{lem:edge-gauge}, the transition between the two sets of block
identifications is specified by one scalar on each upward edge;
the reverse and stationary coefficients are its inverse and $1$.
Thus \eqref{eq:path-gauge} computes the candidate matrices from the
reference matrices. Their bimodule compatibility is already ensured
by transport through the module equivalences.

On a square, denote these transition scalars by $A_s(x),A_s(xr)$
and $B_r(x),B_r(sx)$. The source path from $x$ to $sxr$ first acts
on the right and then on the left, giving $B_r(x)A_s(xr)$; the
target path gives $A_s(x)B_r(sx)$. Since the reference coefficient
is $1$, the candidate coefficient is
\[
h=\frac{B_r(x)A_s(xr)}{A_s(x)B_r(sx)}.
\]
At inputs $sx,xr,sxr$, reversing the appropriate edges gives
$h^{-1},h^{-1},h$. When a path has a stationary step, the source
and target use the same remaining directed edge, so their factors
cancel. If both steps are stationary, both factors are $1$.
This gives the first family.

On a two-point block let the left and right upward edge scalars be $p$ and $q$. Their ratio $k=q/p$ gives the three scalar factors $k^2,k,k$. On the remaining channel, the source path factors are $(k^{-1},1)$ and the target path factors are $(k,1)$; equivalently, applying \eqref{eq:path-gauge} directly gives \eqref{eq:Uk}. Thus every candidate has the asserted form.

Conversely, on a square take $B_r(x)=h$ and all other upward coefficients equal to $1$. On a two-point block take $p=1$, $q=k$. Transporting the reference bimodule through these module equivalences realizes the displayed parameters.
\end{proof}

The lemma describes the interchange in the chosen coordinates.
Now change those coordinates again, denoting the new left and right
upward edge scalars by $A_s(x)$ and $B_r(x)$. The existing parameters
then transform by
\begin{align}
 k&\longmapsto k\frac{B_r(x)}{A_s(x)},\label{eq:gauge-k}\\
 h&\longmapsto h\frac{B_r(x)A_s(xr)}{A_s(x)B_r(sx)}.\label{eq:gauge-h}
\end{align}
These formulas are the coordinate expressions of the module equivalences
in \cref{sec:regular}. A single edge can occur in several joint orbits,
and its edge scalar must be the same in all of them. Thus the
local classification does not by itself choose a simultaneous
normalization; the relations between these local parameters are the
subject of the following sections.

\begin{remark}
The fixed entry $-a$ and the linked coefficients $k^2,k,k$ reflect
compatibility with the full Fibonacci actions. In particular, the
one-loop channels determine $k$, including its sign. The geometric
meaning of this parameter is explained by the digon boundary in
\cref{sec:cochain}.
\end{remark}

\begin{remark}[Inverse interchange matrices]\label{rem:inverse}
The matrix in \eqref{eq:Uk} is invertible for every $k\in\Cx$, and
\begin{equation}\label{eq:Uk-inverse}
U_k^{-1}=\begin{pmatrix}ak^2&bk\\bk&-a\end{pmatrix}.
\end{equation}
This follows directly from \eqref{eq:seed-identities}; $\det U_k=-k^{-2}$. Formula \eqref{eq:Uk-inverse} supplies the inverse for every nonzero interchange parameter $k\in\Cx$. In the unitary specialization $|k|=1$, it also equals the adjoint in the reference metrics.
\end{remark}

\section{Two-generator uniqueness and normalized braid morphisms}\label{sec:rank-two}
We next compare the complete categories generated by two colours.
Their uniqueness will identify the bimodule structures on mixed
three-cells. We also choose an intrinsic normalization of their braid
morphisms, so that pure braid parameters and mixed interchange parameters can later be recorded
in the same cellular complex.

\subsection{The two finite inputs}
For exponents $2$ and $3$, the local parameters of \cref{sec:local}
can be reduced to one scalar $q$. The candidate category already
contains a braid isomorphism between its two alternating simple words.
Its module naturality forces $q=1$. We give this argument here; the
matrices and the two decisive coefficients appear in \cref{app:rank-two}.

\begin{proposition}[Two-generator uniqueness]\label{prop:rank-two}
For $W$ of type $A_1\times A_1$ or $A_2$, every marked algebraic categorification of $\HH(W)$ is marked complex-linearly tensor equivalent to $\Czero(W)$.
\end{proposition}
\begin{proof}
Standardize all complete one-colour actions as in \cref{sec:regular}.
Each joint orbit of an ordered left/right colour pair contributes one
parameter by \cref{lem:local-phases}. For each colour pair, record the number of digons and squares as an
ordered pair. The orbit counts are
\[
\begin{array}{c|cccc|c}
 & (s,s)&(s,t)&(t,s)&(t,t)&\text{parameters}\\\hline
 m=2&(2,0)&(0,1)&(0,1)&(2,0)&6\\
 m=3&(1,1)&(1,1)&(1,1)&(1,1)&8
\end{array}
\]
For $m=2$, left and right multiplication by the same generator coincide,
so each equal-colour pair has two two-point orbits; unequal colours act
transitively on all four vertices. For $m=3$, each colour pair has one
two-point orbit and one four-point orbit in the six-element group.
The corresponding parameters are listed explicitly in
\eqref{eq:rank2-phase-list} and \eqref{eq:rank3-phase-list}.
The exponent matrices for the shared transformations
\eqref{eq:gauge-k}--\eqref{eq:gauge-h} contain square minors of sizes
five and seven with determinants $1$ and $-1$, displayed in
\eqref{eq:M2} and \eqref{eq:M3}. Their inverses have integer entries.
Consequently the chosen edge scalars can be solved for as integer
monomials in the parameters, making the first five or seven parameters
$1$. This normalization requires no root choices and applies to all
nonzero complex parameters, including roots of unity. Write
$q\in\Cx$ for the sole remaining parameter.

Let $A,B$ be the two alternating left words of length $m$.
The corresponding tensor products in the candidate are simple objects
of the same class $h_o$, where $o$ is the longest element. An
isomorphism between these objects induces an invertible module natural
transformation $T:A\to B$ under the regular realization. Rescale it
so that its component at $\eobj{e}$ is $1$ in the fixed simple-object
coordinates. By \cref{lem:evaluation}, every other component is then
forced to be
\begin{equation}\label{eq:forced-braid}
 T_{\eobj{w}}=(\Gamma_{B,w})_{\eobj{e}}^{-1}(\Gamma_{A,w})_{\eobj{e}},
\end{equation}
Here $\Gamma_{A,w}$ and $\Gamma_{B,w}$ are formed from the candidate
interchanges, which still depend on $q$. Thus \eqref{eq:forced-braid}
derives the candidate braid matrix as a function of $q$ from its
normalized value at $\eobj{e}$. All inverses are taken over $\Cx$, using
\eqref{eq:Uk-inverse} for a variable two-dimensional block.

The propagated matrices must satisfy the full module-naturality
identity, since $T$ comes from a morphism in the candidate:
\begin{equation}\label{eq:braid-natural}
 C_{B,r}(TR_r)=(R_rT)C_{A,r}
\end{equation}
for every right colour and every input. Propagation determines the
only possible matrices; this identity tests their consistency with
all right actions. The coefficients in \cref{tab:rank-two-witness}
give, respectively,
\[
 q^2-q=0\quad(m=2),\qquad 1-q=0\quad(m=3).
\]
Since $q\ne0$, both force $q=1$. All elementary interchanges now
equal the reference interchanges, with the complete one-colour
structures fixed. Proposition~\ref{prop:comparison-criterion}
then gives the marked tensor equivalence.
\end{proof}

The role of the calculation is therefore to enforce a necessary
identity of an already existing categorical braid. The exact
coefficients in \cref{app:rank-two} suffice to force $q=1$; at that
value the complete generating actions are those of the reference.

\subsection{An intrinsic normalization of the braid}
A braid isomorphism between simple tensor words is unique only up to
scalar. To compare its action on different parabolic blocks, we need
a normalization preserved by changes of the one-colour module
coordinates. Descent-loop channels provide such a normalization,
because every allowed edge rescaling acts on them as the identity.

Fix $J=\{s,t\}$ with $m=m_{st}\in\{2,3\}$, and let $o_J$ be the longest element of $W_J$. The full tensor subcategory $\cD_J$ generated by $Y_s,Y_t$ has the simple objects $Y_u$, $u\in W_J$: this follows directly from \eqref{eq:hecke-rule} and reduced expressions in $W_J$. Hence \cref{prop:rank-two} applies to $\cD_J$.

The two alternating tensor words of length $m$ are simple and isomorphic. Their intertwiner space is one-dimensional. Choose an invertible braid morphism $B_J$ and apply its left regular action to $Y_{o_J}$. In each of the two length-$m$ paths there is an \emph{all-loop path}: at every step select the descent-loop summand, keeping the intermediate simple equal to $Y_{o_J}$. Normalize $B_J$ by requiring its coefficient between these two paths to be
\begin{equation}\label{eq:braid-loop-normalization}
 \kappa_J=
 \begin{cases}
  1,&m=2,\\
  -a,&m=3.
 \end{cases}
\end{equation}

\begin{lemma}\label{lem:intrinsic-braid}
The normalization \eqref{eq:braid-loop-normalization} exists and determines a unique invertible braid morphism. It is invariant under changes of complete one-colour standardization and is preserved by marked tensor equivalences and regular module transport.
\end{lemma}
\begin{proof}
In the reference model the specified coefficient is $1$ for $m=2$ and the lower-right coefficient $-a$ of $U$ for $m=3$. It is nonzero. A marked algebraic equivalence in \cref{prop:rank-two} carries this normalized morphism to the candidate. Every change of complete one-colour standardization has coefficient $1$ on each descent-loop channel, by \cref{lem:edge-gauge}, and hence leaves the selected coefficient unchanged. Since the intertwiner space is one-dimensional, specifying its nonzero coefficient fixes the morphism. The same reasoning applies to module transport and to equivalences preserving the markings.
\end{proof}

Explicitly, if $\widetilde B_J$ is any nonzero braid morphism and
its all-loop coefficient is $c$, then
\[
 B_J=(\kappa_J/c)\widetilde B_J.
\]
The preceding lemma ensures that $c\ne0$. The normalized morphism
$B_J$ acts coherently on all parabolic cosets by tensoring. Right
braid morphisms are normalized by the right-handed version of the
same prescription, with reversed operation order. The reference
right maps are obtained by reversal and inverse of the real local
blocks; the coefficient \eqref{eq:braid-loop-normalization} is
unchanged.

Here a \emph{pure} polygon uses only left actions or only right actions.
Fix a left parabolic coset $W_Jd$, with $d$ shortest, and write
$A_J,A'_J$ for the source and target tensor words of $B_J$. Their
alternating action lists give two unique ascending paths
$\alpha,\beta$ from $d$ to $o_Jd$. The chosen one-colour module
coordinates determine path identifications
\[
 I_\alpha:Y_{o_Jd}\iso A_J\otimes Y_d,\qquad
 I_\beta:Y_{o_Jd}\iso A'_J\otimes Y_d.
\]
With associativity identifications understood, define
$\theta_{J,d}$ by
\[
 (B_J\otimes\id_{Y_d})I_\alpha=I_\beta\,\theta_{J,d}.
\]
This is a uniquely determined nonzero scalar, since $Y_{o_Jd}$ is
simple and all three maps are invertible. Its value in the reference
coordinates is $1$. The normalization
\eqref{eq:braid-loop-normalization} fixes $B_J$, whereas
$\theta_{J,d}$ records its coefficient in the chosen path
coordinates. For this fixed polygon we abbreviate it to $\theta$,
and orient the boundary by $\alpha$ followed by the reverse of
$\beta$.

For example, when $m=3$ and
$B_J:Y_sY_tY_s\to Y_tY_sY_t$, the two paths are
\[
\begin{aligned}
 \alpha:\quad&
 d\xrightarrow{L_s}sd\xrightarrow{L_t}tsd
   \xrightarrow{L_s}o_Jd,\\
 \beta:\quad&
 d\xrightarrow{L_t}td\xrightarrow{L_s}std
   \xrightarrow{L_t}o_Jd.
\end{aligned}
\]
Under an edge rescaling, write
$I'_\alpha=I_\alpha g(\alpha)$ and
$I'_\beta=I_\beta g(\beta)$, as in \eqref{eq:path-gauge}.
Substituting into the defining equation gives
\[
 \theta'=\theta\,\frac{g(\alpha)}{g(\beta)}.
\]
For the displayed six-sided polygon this is
\[
 \theta'=\theta\,
 \frac{A_s(d)A_t(sd)A_s(tsd)}
      {A_t(d)A_s(td)A_t(std)}.
\]
Thus the change in $\theta$ is the product of edge rescalings
around the oriented polygon. Starting from reference coordinates,
where $\theta=1$, this ratio computes the coefficient in the new
coordinates. The next lemma explains why this single coefficient
also determines the change in every other braid matrix entry on
the same block.

\begin{lemma}[The scalar of a pure braid polygon]\label{lem:polygon-phase}
Fix a coset $W_Jd$ with $d$ shortest and complete one-colour
coordinates, and let $\theta=\theta_{J,d}\in\Cx$ be defined above.
Every nonzero coefficient of the normalized candidate braid equals
the corresponding reference coefficient multiplied by
$\theta^\epsilon$, where $\epsilon\in\{-1,0,1\}$. A nonzero
exponent can occur only when both source and target paths move at
every step. In particular, every coefficient involving a stationary
step is unchanged. Under edge rescalings, $\theta$ is multiplied by
the product of edge scalars around the oriented $2m$-gon. The
right-handed assertion also holds.
\end{lemma}
\begin{proof}
Length additivity identifies $X\mapsto X\otimes Y_d$ with a left
$\cD_J$-module equivalence from the regular two-generator category
onto this coset block: its simple images are precisely the distinct
$Y_{ud}$, $u\in W_J$. Use \cref{prop:rank-two} to identify it with
the reference block. The normalized braid is preserved by
\cref{lem:intrinsic-braid}. Passing from this identification to the
chosen one-colour coordinates therefore changes its matrix only by
the edge rescalings of \cref{lem:edge-gauge}. By
\eqref{eq:path-gauge}, each coefficient is multiplied by the source
path product divided by the target path product. On the two
ascending paths, this ratio is precisely $\theta$.

For any matrix entry, join its source path to the reverse of its
target path. After discarding stationary steps, this gives a closed
walk of length at most $2m$ in a single $2m$-gon. Its signed edge
chain is an integer multiple of the oriented polygon boundary.
A nonzero winding requires at least $2m$ moving steps, so that
multiple is $\epsilon\in\{-1,0,1\}$. If $\epsilon\ne0$, both
length-$m$ paths move at every step and follow opposite
half-boundaries; the multiplier is then $\theta^\epsilon$.
Otherwise the signed edge chain is zero and the coefficient is
unchanged. This also covers the coefficients in the two-dimensional
block. A further edge rescaling changes $\theta$ by the boundary
product derived above. The normalization
\eqref{eq:braid-loop-normalization} fixes the overall scalar.
\end{proof}

\section{The diagonal quotient of two Davis complexes}\label{sec:geometry}
On an individual block, a reference identification makes its face
parameters equal to $1$. The remaining problem is to choose such
identifications simultaneously: an edge may belong to several
faces, and its rescaling must be the same in each occurrence.
We organize this problem in the diagonal quotient of
two Davis complexes. Its vertices are the simple labels, its edges
carry the left and right edge rescalings, and its two-cells carry
$h,k,\theta$. Mixed three-cells will encode the categorical
compatibility among these parameters.

We establish the geometry before assigning the scalars: quotient
characteristic domains must be balls in order to use cellular chains,
and contractibility will provide the fillings used in the simultaneous
normalization. Throughout, the chains are those of the ordinary quotient,
with its actual integer incidence coefficients.

\subsection{Product-cell orbits and the quotient}
Use the Davis complex $\Sigma$ and its Coxeter cells $wP_J$ from
\cref{subsec:davis-prelim}. Define
\begin{equation}\label{eq:Y}
 Y=(\Sigma\times\Sigma)/W,
\end{equation}
where $W$ acts diagonally on the left. At vertices, the identification
\begin{equation}\label{eq:vertex-id}
 [(u,v)]\longmapsto u^{-1}v
\end{equation}
identifies the vertex set of $Y$ with $W$. A first-factor edge $u\to us$ induces $w\to sw$, and a second-factor edge $v\to vr$ induces $w\to wr$. We retain the left/right distinction even when these edges have the same endpoints. Thus $Y^{(1)}$ is exactly the graph on which the shared edge rescalings are defined.

An orbit of product cells has a representative $P_J\times dP_K$, where $J,K$ are spherical and $d$ is the shortest representative of $W_JdW_K$. Its setwise stabilizer is
\begin{equation}\label{eq:H}
 H=W_J\cap dW_Kd^{-1}.
\end{equation}
We denote the corresponding quotient cell by
\[
 \cell{J}{K}{d},\qquad\dim\cell{J}{K}{d}=|J|+|K|.
\]
A characteristic domain is $(P_J\times dP_K)/H$. The possibility $H\ne\{e\}$ must be treated before using ordinary cellular chains.

\begin{lemma}[Quotient characteristic balls]\label{lem:quotient-balls}
The pair
\[
 \bigl((P_J\times dP_K)/H,\,\partial(P_J\times dP_K)/H\bigr)
\]
is homeomorphic to $(D^n,\partial D^n)$, where $n=|J|+|K|$.
For $n=0$, this means the pair $(\{\mathrm{pt}\},\varnothing)$.
Product-cell orbits define a locally finite CW structure on $Y$.
\end{lemma}
\begin{proof}
When $n=0$, both factors are points and the ball-pair assertion is
immediate. For $n>0$, we identify the stabilizer representation and
construct a homeomorphism of the quotient ball pair. We then verify
the cellular structure in all dimensions.

Use the Euclidean models of the two cells from \cref{subsec:davis-prelim}.
Identify the copy $dP_K$ with $P_K$ by $d\cdot y\mapsto y$. Since
$h(d\cdot y)=d\cdot((d^{-1}hd)y)$, an element $h\in H$ acts in the
two sets of local coordinates by
\[
 (x,y)\longmapsto\bigl(hx,(d^{-1}hd)y\bigr).
\]
In each factor the action permutes the vertices and hence fixes their
barycentre. Taking these barycentres as origins gives linear actions
on the two spaces of displacement vectors.

Intersections of parabolic subgroups are parabolic; the intersection is parabolic in each of the two groups in \eqref{eq:H}. This follows from the parabolic intersection theorem and its application inside the corresponding finite reflection groups; see \cite[Lemma 5.3.6]{Davis}. In each Euclidean reflection representation, restriction to $H$ splits into the reflection representation $E$ of $H$ and a trivial orthogonal complement. The two copies of $E$ agree as $H$-representations: both identify with the reflection subspace of $H$ in the representation $V$ of \cref{subsec:geometric-representation}. The resulting representation on the product of the two displacement spaces is therefore
\begin{equation}\label{eq:doubled}
 E\oplus E\oplus\RR^f\cong E_{\CC}\oplus\RR^f.
\end{equation}
The trivial case $H=\{e\}$ is included.

The action of $H$ on $E_{\CC}$ is generated by complex reflections. By the Chevalley--Shephard--Todd theorem \cite{Chevalley,ShephardTodd}, its polynomial invariant algebra has homogeneous algebraically independent generators $f_1,\ldots,f_k$, where $k=\dim_{\RR}E$. These define the topological orbit-space homeomorphism $E_{\CC}/H\cong\CC^k$: invariant polynomials separate finite orbits, and the finite algebraic quotient map is proper. Only the classical characteristic-zero form of the theorem is used; see also \cite[Theorem 1.1]{Satriano}.

To control the boundary sphere as well as the ambient orbit space, put $d_i=\deg f_i$ and define coordinate homeomorphisms
\[
 \chi_i(z)=
 \begin{cases}
 |z|^{1/d_i-1}z,&z\ne0,\\
 0,&z=0.
 \end{cases}
\]
The map induced by $(\chi_1(f_1),\ldots,\chi_k(f_k))$ is positively homogeneous of degree one for real positive dilation. Extend it by the identity on $\RR^f$. After radial normalization it identifies the quotient unit sphere in \eqref{eq:doubled} with the ordinary unit sphere. Coning this sphere homeomorphism identifies the quotient ball and its boundary. The product $P_J\times dP_K$, centred at its $H$-fixed centre, is equivariantly radially homeomorphic to such a ball, giving the required pair.

It remains to check that these ball pairs are characteristic domains
for the quotient cells. The action on the product complex is cellular
and proper, and the product complex is locally finite because $S$ is finite. Each quotient interior is embedded in $Y$: if a translate of an open product cell meets it, the translate is that same cell, so the element belongs to its setwise stabilizer. An open product-cell orbit is the interior of the quotient ball just established. Each closure is incident to finitely many lower-dimensional cell orbits, and the quotient topology is the weak topology for these characteristic maps. They therefore define the stated CW structure. The attaching maps retain the boundary-face identifications induced by the quotient.
\end{proof}

\begin{lemma}[Contractibility]\label{lem:contractible}
The ordinary quotient $Y$ is contractible.
\end{lemma}
\begin{proof}
Use the $W$-invariant CAT(0) metric recalled in
\cref{subsec:davis-prelim}. Let $\gamma_{y,x}:[0,1]\to\Sigma$ be the unique geodesic from $y$ to $x$. Geodesics depend continuously and $W$-equivariantly on their endpoints, so
\[
 H_t([x,y])=[x,\gamma_{y,x}(t)]
\]
is a strong deformation retraction onto the image of the diagonal. That image is the Davis chamber $\Sigma/W\cong\mathcal K$, which is
contractible by \cref{subsec:davis-prelim}.
\end{proof}

Contractibility will be used in two degrees: every finite one-cycle
bounds a finite two-chain, and every finite two-cycle bounds a finite
three-chain. We use the cellular chain groups of the CW space $Y$:
\[
 C_k(Y;\ZZ)=\bigoplus_{\substack{F\text{ an oriented cell of }Y\\\dim F=k}}\ZZ[F].
\]
Each chain is a finite integer sum of quotient cells. The coefficient of
a $(k-1)$-cell in $\partial[F]$ is the degree of the attaching map on
that face, computed in the quotient characteristic ball. This convention
includes cells with nontrivial stabilizers and determines the integral
incidences used below.

\section{The scalar cochain and mixed three-cell compatibility}\label{sec:cochain}
We now express the local parameters and their rescalings in the cellular coordinates
of $Y$. The local parameters form a multiplicative two-cochain, and
changing module coordinates multiplies it by a coboundary. To obtain
compatibility on a mixed three-cell, we compare the entire bimodule
on that cell with the reference; one such bimodule equivalence controls all
its boundary faces at once.

\subsection{The two-dimensional dictionary}
Orient every left and right coloured edge upward in length and denote
these oriented cells by $\mathsf e_s^L(x)=[x\to sx]_L$ and
$\mathsf e_r^R(x)=[x\to xr]_R$. For an edge cochain
$u:C_1(Y;\ZZ)\to\Cx$, their values are the edge scalars:
\[
 A_s(x)=u(\mathsf e_s^L(x)),\qquad
 B_r(x)=u(\mathsf e_r^R(x)).
\]
Thus edges are added in cellular boundaries, while their scalar values
are multiplied under $u$. A cell is called \emph{mixed} when both $J$ and $K$ are nonempty, and \emph{pure} otherwise. There are three types of two-dimensional cells.

If $|J|=|K|=1$ and the stabilizer is trivial, the cell is a mixed square. For $J=\{s\}$, $K=\{r\}$ and shortest vertex $x$, orient it by the path which first moves right and then left, minus the path which first moves left and then right. In edge-chain notation its boundary is
\begin{equation}\label{eq:square-boundary}
 \mathsf e_r^R(x)+\mathsf e_s^L(xr)-\mathsf e_s^L(x)-\mathsf e_r^R(sx).
\end{equation}
Its scalar is the $h$ in \cref{lem:local-phases}.

If the stabilizer is nontrivial, it has order two and exchanges the opposite vertices of the product square. The characteristic domain is a square modulo a half-turn, hence a disk. Its boundary is a digon with the two differently coloured edges from $x$ to $y=sx=xr$. Choose the orientation
\begin{equation}\label{eq:digon-boundary}
 \mathsf e_r^R(x)-\mathsf e_s^L(x).
\end{equation}
Its scalar is $k$. The all-moving scalar channel in \cref{lem:local-phases}(ii)
follows this quotient boundary twice and has coefficient $k^2$.
A channel with one descent loop follows it once and has coefficient
$k$. Thus the face value is $k$: the full interchange detects the
primitive boundary, including the sign information that its square
alone would lose.

Finally, cells of types $(2,0)$ and $(0,2)$ are pure braid polygons. Orient their boundaries by the source braid path followed by the reverse of the target path, and assign the scalar $\theta$ of \cref{lem:polygon-phase}.

\begin{figure}[ht]
\centering
\begin{tikzpicture}[>=Stealth,scale=0.93,every node/.style={font=\small}]
 \node (x) at (0,0) {$x$};
 \node (sx) at (0,1.8) {$sx$};
 \node (xr) at (2.8,0) {$xr$};
 \node (sxr) at (2.8,1.8) {$sxr$};
 \draw[->] (x)--node[left] {$\mathsf e_s^L(x)$}(sx);
 \draw[->] (xr)--node[right] {$\mathsf e_s^L(xr)$}(sxr);
 \draw[->] (x)--node[below] {$\mathsf e_r^R(x)$}(xr);
 \draw[->] (sx)--node[above] {$\mathsf e_r^R(sx)$}(sxr);
 \node at (1.4,.9) {$h$};
 \node at (1.4,-.8) {trivial stabilizer};
 \node (u) at (6.1,0) {$x$};
 \node (v) at (6.1,1.8) {$y$};
 \draw[->] (u) to[bend left=48] node[left] {$\mathsf e_s^L(x)$} (v);
 \draw[->] (u) to[bend right=48] node[right] {$\mathsf e_r^R(x)$} (v);
 \node at (6.1,.9) {$k$};
 \node at (6.1,-.8) {order-two stabilizer};
\end{tikzpicture}
\caption{Mixed quotient faces. The primitive scalar on the digon is $k$, even though one matrix channel has coefficient $k^2$. All arrows are oriented upward in Coxeter length.}
\label{fig:faces}
\end{figure}

Assigning the inverse to the opposite orientation and extending by
$c(\sum n_FF)=\prod c(F)^{n_F}$, these scalars define a multiplicative
cellular cochain
\begin{equation}\label{eq:phase-cochain}
 c:C_2(Y;\ZZ)\longrightarrow\Cx.
\end{equation}
A family of edge rescalings is exactly a cochain $u:C_1(Y;\ZZ)\to\Cx$. Equations \eqref{eq:gauge-k}, \eqref{eq:gauge-h} and \cref{lem:polygon-phase} give
\begin{equation}\label{eq:cochain-gauge}
 c\longmapsto c\,\delta u,\qquad (\delta u)(F)=u(\partial F).
\end{equation}
For example, applying $u$ to \eqref{eq:square-boundary} gives
\[
 u(\partial F)=\frac{B_r(x)A_s(xr)}{A_s(x)B_r(sx)},
\]
which is precisely the change of the square matrix parameter in
\eqref{eq:gauge-h}.

\subsection{Compatibility on mixed three-cells}
A mixed three-cell has type $(2,1)$ or $(1,2)$: it combines a finite
two-generator action with a commuting one-generator action. We use
\cref{prop:rank-two} to identify this full bimodule, including the
normalized braid on its pure faces. The resulting edge rescalings
then account for every boundary face by the same coboundary formula.

\begin{lemma}[Mixed cocycle equations]\label{lem:mixed-cocycle}
For every mixed three-cell $Q$,
\begin{equation}\label{eq:mixed-closed}
 c(\partial Q)=1.
\end{equation}
\end{lemma}
\begin{proof}
Consider $Q=\cell{J}{\{r\}}{d}$ with $J=\{s,t\}$ spherical and $d$ shortest in $W_JdW_r$. Apply \cref{lem:transport} to $d,r$.

Suppose first that $dr$ is shortest in its left $W_J$-coset. The two cosets are distinct, and right minimality of $d$ gives $\ell(dr)=\ell(d)+1$. Their union is the whole double block. For $u\in W_J$ and $\epsilon\in\{0,1\}$, the products
\[
 Y_u\otimes Y_d\otimes Y_r^{\otimes\epsilon}
\]
are distinct simple objects, labelled by $ud$ and $udr$. Indeed \eqref{eq:length-add} gives length additivity in both cosets. Consequently
\[
 (A,B)\longmapsto A\otimes Y_d\otimes B
\]
is an equivalence of the complete $\cD_J$--$\Fib$ bimodule on the double block with the regular independent-factor model. A marked algebraic equivalence from \cref{prop:rank-two} identifies $\cD_J$ with the reference two-generator category.

In the second case, $dr=jd$ for a unique $j\in J$, and the block
is $W_Jd$. The functor $X\mapsto X\otimes Y_d$ identifies its left
$\cD_J$-action with the regular left module. For a tensor category
$\mathcal{A}$, endofunctors of its regular left module are right tensoring
functors. Their composition reverses the order of tensor factors, so
\[
 \operatorname{End}_{\mathcal{A}}(\mathcal{A})\simeq\mathcal{A}^{\mathrm{rev}},
\]
where $\mathcal{A}^{\mathrm{rev}}$ has the reversed tensor product.
A commuting right $\Fib$-action is likewise a tensor functor from
$\Fib^{\mathrm{rev}}$ to this endofunctor category. With
$\mathcal{A}=\cD_J$, reversing both tensor products therefore gives a tensor
functor $\Fib\to\cD_J$.

Evaluation at the regular unit, which corresponds to $Y_d$ on the
block, sends $\tau$ to $Y_j$, since
$Y_d\otimes Y_r\cong Y_j\otimes Y_d$. The functor lands in the
marked subcategory $\Fib_j$. Comparing it with the marked inclusion
reduces to a tensor autoequivalence of $\Fib$, which is monoidally
trivial by \cref{lem:fib-aut}. Together with the two-generator
equivalence of \cref{prop:rank-two}, this identifies the entire
bimodule with the reference one.

In either case, choose this single bimodule equivalence on the
whole block and identify its underlying labelled linear functor with
the identity. Its restrictions to all one-colour blocks have the
rescaling form of \cref{lem:edge-gauge}. They assign one scalar to each
coloured edge of the block, even if that edge belongs to several
boundary faces. Denote by $v$ the local edge cochain for the transition
from reference coordinates to candidate coordinates.

On every mixed two-face, \eqref{eq:gauge-k}--\eqref{eq:gauge-h}
express its scalar as $v$ evaluated on the boundary. On every pure
two-face the same assertion follows from \cref{lem:polygon-phase}:
the bimodule equivalence preserves the normalized braid by
\cref{lem:intrinsic-braid}. Thus $c(F)=v(\partial F)$ on all
boundary faces of $Q$, with the same $v$. Using the actual integer
quotient incidences gives
\[
 c(\partial Q)=v(\partial\partial Q)=1.
\]
Repeated or identified faces are counted with their incidence
coefficients, so the equation also covers nontrivial stabilizers.
The type $(1,2)$ case interchanges left and right; regular right-module
endofunctors recover the tensor product in its original order.

\end{proof}

The categorical argument has established \eqref{eq:mixed-closed}
on mixed three-cells. To obtain triviality on every integral two-cycle,
we must also account for pure three-cell boundaries. The next section
expresses them as integral combinations of mixed boundaries.

\section{Integral elimination of pure cellular boundaries}\label{sec:integral}
The categorical argument gives $c(\partial Q)=1$ for mixed
three-cells. Contractibility of $Y$ expresses every two-cycle as a
sum of three-cell boundaries, but these may include pure cells.
We bridge this gap by replacing every pure boundary with an integer
combination of mixed boundaries. The unit incidence coefficients in
this replacement are what allow the later normalization to use
integer powers of the parameters, without extracting roots.

Let $C_k^{\mix}(Y;\ZZ)$ be the subgroup generated by mixed $k$-cells.
We restrict the cellular boundary of $C_k(Y;\ZZ)$ to this subgroup;
its image may contain pure faces.

\begin{theorem}[Pure-boundary elimination]\label{thm:pure-elimination}
For every $k\ge2$,
\begin{equation}\label{eq:pure-elimination}
 \partial C_k(Y;\ZZ)=\partial C_k^{\mix}(Y;\ZZ).
\end{equation}
More precisely, the boundary of each pure $k$-cell is a finite integral combination of boundaries of mixed $k$-cells. Every elimination of a pure facet uses a coefficient $1$ or $-1$.
\end{theorem}
\begin{proof}
It suffices to consider a pure left cell
\[
 P_d=\cell{J}{\varnothing}{d},\qquad |J|=k,
\]
where $J$ is spherical and $d$ is shortest in $W_Jd$. The right case
reverses the factors. We first reduce the length of $d$ by ordinary
product prisms, and then eliminate the remaining boundary at $d=e$
using one quotient prism.

Assume $d\ne e$. Choose a right descent $r$ of $d$ and put $d'=dr$, so $\ell(d')=\ell(d)-1$. Such an $r$ exists as the last letter of a reduced expression for $d$. The alternative $dr=jd$ with $j\in J$ in \cref{lem:transport} is impossible: $d$ has no left $J$-descent, so it would give $\ell(dr)=\ell(jd)=\ell(d)+1$. Hence $d'$ is shortest in its left coset and $W_Jd'$ is different from $W_Jd$. Since $d'r=d$ is longer, $d'$ is also the shortest double-coset representative for the following cell:
\[
 Q=\cell{J}{\{r\}}{d'}.
\]
Its stabilizer is trivial; a nontrivial stabilizer would identify those two left cosets. Thus $Q$ is a product prism $P_J\times I$.

The two end facets are precisely $P_d$ and $P_{d'}$, and each has
incidence coefficient $1$ or $-1$. The product boundary formula is
\[
 \partial(P_J\times I)
 =(\partial P_J)\times I+(-1)^kP_J\times(\{1\}-\{0\}).
\]
Every side facet is $F\times I$ for a $(k-1)$-face $F$ of $P_J$.
It therefore has type $(k-1,1)$, with both factors positive-dimensional
because $k\ge2$. These are exactly the mixed facets; see
\cref{fig:elimination-prism}.

\begin{figure}[ht]
\centering
\begin{tikzpicture}[>=Stealth,scale=.95,every node/.style={font=\small}]
 \coordinate (a) at (0,0); \coordinate (b) at (3.0,0);
 \coordinate (c) at (4.0,.8); \coordinate (d) at (1.0,.8);
 \coordinate (A) at (0,2.3); \coordinate (B) at (3.0,2.3);
 \coordinate (C) at (4.0,3.1); \coordinate (D) at (1.0,3.1);
 \fill[blue!9] (A)--(B)--(C)--(D)--cycle;
 \fill[orange!12] (a)--(b)--(B)--(A)--cycle;
 \draw[dashed] (a)--(d)--(c); \draw[dashed] (d)--(D);
 \draw (a)--(b)--(c)--(C)--(D)--(A)--(a);
 \draw (A)--(B)--(C); \draw (b)--(B);
 \node at (2,2.7) {$P_d\quad(k,0)$};
 \node at (1.6,-.45) {$P_{d'}\quad(k,0)$};
 \node[align=center] at (1.5,1.25) {$F\times I$\\mixed: $(k-1,1)$};
 \draw[->] (4.8,.4)--node[right] {$I$}(4.8,2.7);
\end{tikzpicture}
\caption{The product prism in the length-reduction step, illustrated for
$k=2$. The two horizontal facets are pure. Each side facet is a face of
$P_J$ crossed with $I$, hence is mixed. For $k=3$, the same construction
is a four-dimensional prism with two pure three-dimensional ends and
mixed facets of type $(2,1)$.}
\label{fig:elimination-prism}
\end{figure} Write
\begin{equation}\label{eq:prism-boundary}
 \partial Q=\epsilon P_d+\epsilon'P_{d'}+M,
 \qquad \epsilon,\epsilon'\in\{1,-1\},\quad M\in C_k^{\mix}(Y;\ZZ).
\end{equation}
Taking the boundary and solving using the unit $\epsilon$ gives
\begin{equation}\label{eq:length-recursion}
 \partial P_d=-\epsilon\epsilon'\partial P_{d'}-\epsilon\partial M.
\end{equation}
Repeating this step strictly decreases $\ell(d)$, and therefore reduces the assertion in finitely many steps to $d=e$.

For the last step choose $r\in J$ and consider
\begin{equation}\label{eq:closing-prism}
 Q_0=\cell{J}{\{r\}}{e}
       =(P_J\times I)/\langle(r,r)\rangle.
\end{equation}
Here the first $r$ is a reflection of $P_J$ and the second reverses $I$. By \cref{lem:quotient-balls} the characteristic domain is a ball. The involution exchanges the two end facets. Their common image is one pure $k$-facet $P_e$, with incidence coefficient $1$ or $-1$.

For clarity, choose one end facet. Its interior has trivial stabilizer and maps homeomorphically to the interior of $P_e$; it represents that quotient facet once. The two reflections together preserve the ambient orientation, so the contributions from the exchanged end facets have the same induced sign. Equivalently, the ambient quotient has degree two and the two lifts of the pure facet have total degree two, leaving unit incidence in the quotient ball. All other $k$-facets are images of side facets and are mixed. Thus
\begin{equation}\label{eq:closing-boundary}
 \partial Q_0=\epsilon_0P_e+M_0,
 \qquad\epsilon_0\in\{1,-1\},\quad M_0\in C_k^{\mix}(Y;\ZZ).
\end{equation}
Taking a second boundary yields
$\partial P_e=-\epsilon_0\partial M_0$. At this last step, as in
\eqref{eq:length-recursion}, the coefficient of the boundary being
eliminated is a unit in $\ZZ$. Combining the last step with the finite
length recursion proves \eqref{eq:pure-elimination} over $\ZZ$.
\end{proof}

\begin{remark}\label{rem:no-larger-parabolic}
Each step uses a spherical $J$-cell in one factor and a rank-one cell
in the other. These two finite cells exist for every choice of the
external colour $r$, so the length recursion also applies when $W$ is
infinite. In the final quotient prism, the end facets form a single
orbit and its quotient incidence is a unit. This is the reason that the
elimination works over $\ZZ$.
\end{remark}

\begin{corollary}[Integral generation of two-cycles]\label{cor:two-cycles}
\begin{equation}\label{eq:two-cycles}
 \ker\bigl(\partial:C_2(Y;\ZZ)\to C_1(Y;\ZZ)\bigr)
 =\partial C_3^{\mix}(Y;\ZZ).
\end{equation}
\end{corollary}
\begin{proof}
By \cref{lem:contractible}, $\ker\partial_2=\im\partial_3$. Apply \cref{thm:pure-elimination} with $k=3$ to replace every pure three-cell boundary by mixed three-cell boundaries.
\end{proof}

In particular, the boundary of a pure three-colour cell is an integral
combination of mixed boundaries. The two-generator marked equivalences
of \cref{lem:mixed-cocycle} therefore supply the compatibility on every
three-cell. This is the passage from local two-generator data to the global
normalization.

\section{Simultaneous normalization and marked tensor equivalence}\label{sec:global}
The preceding sections show that $c$ takes the value $1$ on every
integral two-cycle. We now construct a single family of edge rescalings that removes all
face parameters. We set the scalars on a spanning tree to $1$;
filling the remaining edge cycles determines the other scalars. We then use the regular realization
to recover a marked tensor equivalence of the original categories.

\subsection{Simultaneous edge rescalings}
\begin{theorem}[Simultaneous normalization]\label{thm:normalization}
For the complete left and right actions arising from any marked candidate $\cD$, there is a family of edge rescalings which sends every mixed interchange to the reference blocks $1,U$ and every normalized pure braid scalar to $1$. For finite or infinite $W$, the resulting module-functor constraints are invertible linear maps on every finite-support object.
\end{theorem}
\begin{proof}
By \cref{lem:mixed-cocycle,cor:two-cycles}, the cochain $c$ is trivial
on every integral two-cycle. We construct a one-cochain $u$ with
$\delta u=c$, verify that it is independent of the chosen fillings,
and interpret its inverse as a family of edge rescalings.

Fix a total order on the finite set $S$. For each $w\ne e$, choose its
lexicographically least reduced expression $w=s_1\cdots s_k$, put
$p(w)=s_1\cdots s_{k-1}$, and retain the right edge
\[
 p(w)\xrightarrow{\ s_k\ }w.
\]
The choice exists because the set of reduced expressions of $w$ is finite.
Every prefix is again lexicographically least for the element it
represents: replacing it by a smaller reduced expression would give a
smaller reduced expression of $w$ of the same length. These parent edges
form a spanning tree $T$ rooted at $e$. Indeed, iterating $p$ decreases
length by one until it reaches $e$, and a vertex of maximal length on a
cycle would have two distinct parent edges. The tree path from $e$ to
$w$ is its chosen reduced expression and has length $\ell(w)$.
This defines $T$ vertex by vertex for finite or infinite $W$.

The edges of $Y^{(1)}$ retain their left or right type and their colour,
even when two edges have the same endpoints. Thus $T$ uses only right edges;
all left edges and the remaining right edges lie outside $T$.
For each oriented edge $e_0$ outside $T$, let $z_{e_0}$ be the finite cycle consisting of $e_0$ followed by the tree path from its endpoint back to its starting point. Since $H_1(Y;\ZZ)=0$, choose a finite integer two-chain $F_{e_0}$ with
\[
 \partial F_{e_0}=z_{e_0}.
\]
Define
\begin{equation}\label{eq:global-u}
 u(e_0)=c(F_{e_0}),\qquad u(e)=1\text{ for tree edges},
\end{equation}
and use inverses for reverse orientations. If two fillings are chosen, their difference is a finite two-cycle, on which $c$ is $1$. Thus \eqref{eq:global-u} is independent of the filling.

Fundamental tree cycles form a basis of the group of finite graph one-cycles. For any two-cell $F$, write
\[
 \partial F=\sum_{e_0\notin T}n_{e_0}z_{e_0},
\]
with finitely many nonzero integers $n_{e_0}$. The chain
$F-\sum n_{e_0}F_{e_0}$ is a two-cycle. Therefore
\[
 c(F)=\prod_{e_0\notin T}c(F_{e_0})^{n_{e_0}}
      =u(\partial F),
\]
so $\delta u=c$.

To return to the categorical data, assign $u(\mathsf e_s^L(x))$ and
$u(\mathsf e_r^R(x))$ to the corresponding ascending edges. By
\cref{lem:edge-gauge}, every such scalar determines a complete
module-functor structure on its regular Fibonacci block. Taking the
sum over the blocks supplies a module autoequivalence with underlying
functor $\Id_{\cM}$ for each left and right colour. We apply the inverse family: by \eqref{eq:cochain-gauge},
the new cochain is $c(\delta u)^{-1}=1$. The local classification
and \cref{lem:polygon-phase} then identify all mixed blocks and
normalized pure braid coefficients with those of the reference.

The construction applies to finite and infinite $W$ alike. Each $u(e)$
is a finite product of the original scalars with integer exponents,
because its filling chain and tree path are finite. Every object of
$\cM$ has finite support, so these scalars define invertible
module-functor constraints on each finite-support object.
\end{proof}

\subsection{Recovering the full tensorator}
\label{subsec:global-tensorator}
\begin{proof}[Proof of \cref{thm:main}]
Apply \cref{thm:normalization}. Write $\beta_{r,Z}$ for the resulting right coordinate changes, and let $Q\in\cE_0$ have underlying functor $\Id_{\cM}$ and constraints $\beta_{r,Z}^{-1}$. Set
\[
P=\operatorname{Ad}_Q\circ\operatorname{Reg}_{\cD},\qquad R=\operatorname{Reg}_0.
\]
By \cref{prop:comparison-criterion}, the normalized left block identifications
identify the marked seed functors in these realizations, and $P,R$
have the same replete full tensor image $\mathcal H\subseteq\cE_0$.

We spell out how this determines the tensorator. Regard
\[
P:\cD\longrightarrow\mathcal H,\qquad
R:\Czero(W)\longrightarrow\mathcal H
\]
as fully faithful strong monoidal equivalences. Choose a linear equivalence $E:\Czero(W)\to\cD$ and a linear natural isomorphism
\[
\sigma:PE\iso R.
\]
Such a choice can be made for any of the label-preserving linear equivalences from \cref{subsec:regular-target}: both sides send each simple label to the same simple isomorphism class in $\mathcal H$.

Write the tensorators of $P$ and $R$ as
\[
p_{X,Y}:P(X)P(Y)\iso P(X\otimes Y),\qquad
r_{A,B}:R(A)R(B)\iso R(A\otimes B).
\]
For this choice of $\sigma$, full faithfulness of $P$ defines a unique map $J_{A,B}$ by
\begin{equation}\label{eq:global-tensorator}
P(J_{A,B})=
\sigma_{A\otimes B}^{-1}\,
 r_{A,B}(\sigma_A\otimes\sigma_B)\,
 p_{E(A),E(B)}^{-1}.
\end{equation}
Its source and target are $P(E(A)\otimes E(B))$ and $P(E(A\otimes B))$, respectively. The right side is an isomorphism, so $J_{A,B}$ is an isomorphism. If $p_0$ and $r_0$ are the unit maps of $P$ and $R$, the unit isomorphism $J_0:\mathbf1_{\cD}\to E(\mathbf1)$ is defined by
\[
P(J_0)=\sigma_{\mathbf1}^{-1}r_0p_0^{-1}.
\]
Naturality follows from that of $\sigma,p,r$. Applying $P$ to the two sides of \eqref{eq:tensorator-coherence}, and using \eqref{eq:global-tensorator}, reduces the equation to the coherence of $P$ and $R$. Faithfulness of $P$ then gives the equation in $\cD$. The unit axioms follow in the same way. Thus $(E,J,J_0)$ is a strong monoidal equivalence and $\sigma$ is monoidal.

The normalized seed identifications give monoidal natural isomorphisms $R\iota_s^0\cong P\iota_s$. Composing with $\sigma\iota_s^0$ and lifting through the fully faithful $P$ gives monoidal natural isomorphisms
\[
E\iota_s^0\cong\iota_s.
\]
The equivalence is therefore marked. The simple labels are preserved by construction, or equivalently by the reduced tensor words in the marked generators.
\end{proof}

The local edge rescalings determine identifications of the generating module actions; \eqref{eq:global-tensorator} supplies their coherent extension to arbitrary tensor products. Throughout the argument, the tensor product and associator of $\cD$ remain fixed. The normalization changes the chosen module identifications, and the common regular image converts those choices into the required marked tensor equivalence.

\section{ADE classification and further consequences}\label{sec:consequences}
The rigidity, construction, and uniqueness results now give the
classification announced in the Introduction. We first explain its
finite ADE specialization, then discuss unitary structures, seed
branches, and infinite type. We conclude with an integral invariant
of the simultaneous normalization.

\subsection{The finite ADE classification}
For a category with marked positive Fibonacci generators,
\cref{thm:ring-rigidity,cor:exponent-restriction} recover a Coxeter
system satisfying \eqref{eq:exponents}. Construction and marked
uniqueness, \cref{thm:construction,thm:main}, identify the category
with $\Czero(W)$. Thus the permitted matrices classify the marked
tensor categories, as stated in \cref{cor:classification}.

The simple objects of $\Czero(W)$ are indexed by $W$, so finite
categorical rank is equivalent to finiteness of this group. Among
the allowed systems, the classical criterion in
\cref{subsec:geometric-representation} selects precisely the disjoint
unions of finite ADE diagrams. This is the finite classification of
\cref{cor:ADE}: the diagram specifies the marked generators, the group
elements label the simples, and Coxeter lengths determine the dimensions
in \eqref{eq:dimension-main}. The classical Coxeter classification supplies
the list of finite diagrams; the construction and uniqueness theorems
supply their categorical realizations and exhaustiveness.

\subsection{Unitarizability and unitary tensor equivalence}
There are two related consequences for unitary structures. First,
the tensor equivalence transports the reference $C^*$-structure to
an algebraic candidate. Second, when unitary structures and unitary
markings are already given, the block identifications and subsequent rescalings can be chosen to
respect them.

\begin{corollary}[Unitarizability]\label{cor:unitarizable}
For finite $W$, every marked categorification in \cref{def:marked} is unitarizable.
\end{corollary}
\begin{proof}
Transport the $C^*$-tensor structure of the reference along the tensor equivalence in \cref{thm:main}. One may first replace the equivalence by an identification of labelled semisimple skeletons, and then pull back the adjoint, norms and tensor constraints. The transported structure is a compatible $C^*$-tensor structure on the candidate.
\end{proof}

\begin{corollary}[Marked unitary tensor equivalence]\label{cor:unitary-comparison}
Suppose that candidate and reference are equipped with unitary structures and that their seed embeddings are unitary. Then the equivalence in \cref{thm:main} can be chosen unitary, with unitary monoidal natural isomorphisms on the markings. The elementary block rescalings have modulus one, also when $W$ is infinite.
\end{corollary}
\begin{proof}
Standardize one-colour actions by unitary module equivalences. The initial edge and cell parameters lie in the unit circle $\TT=\{z\in\CC:|z|=1\}$. The two-generator elimination uses unimodular integer exponent matrices, so its integer monomials preserve $\TT$. The spanning-tree primitive in \eqref{eq:global-u} is also a finite integer monomial in the initial parameters. Thus the resulting block identifications are unitary module equivalences, and the regular images agree as generated $C^*$-tensor categories. Choose the underlying linear equivalences to preserve adjoints,
and choose the simple-object identifications in $\sigma:PE\to R$
to be unitary, extending them over finite direct sums. The tensorators
of the regular realizations and the conjugation are then unitary, so \eqref{eq:global-tensorator} and its unit analogue lift to
unitary maps. The monoidal natural isomorphisms on the seed functors lift in the same way.
This gives the claimed marked unitary equivalence. Every elementary
edge scalar has modulus one, including in infinite rank.
\end{proof}

In finite rank, the unitary tensor equivalence also follows from Reutter's general theorem for unitarizable categories \cite{Reutter}. The explicit normalization above gives such an equivalence for infinite $W$ as well.

\subsection{Tensor branches and connected components}\label{subsec:branches}
The fusion rule $\tau^2=\one+\tau$ has two tensor realizations: the
positive Fibonacci category and the Yang--Lee category
\cite[Section~3]{BookerDavydov}. Put
\[
 d_+=\ph,\qquad d_-=-\ph^{-1},\qquad
 U_d=\begin{pmatrix}d^{-1}&b_d\\b_d&-d^{-1}\end{pmatrix},
 \quad b_d^2=d^{-1}.
\]
These are the two Galois-conjugate seeds. The Yang--Lee category is
nonunitarizable. Their intrinsic squared norms at the nontrivial simple
are $d_+^2=\ph^2$ and $d_-^2=\ph^{-2}$, respectively.

The braid relation forces the branches to agree across every exponent-three
edge. To see this under our algebraic hypotheses, write $n(X)$ for the
intrinsic squared norm of a simple object. It is nonzero, and
\[
 n(XY)=n(X)n(Y)\qquad\text{if }X,Y,XY\text{ are simple}.
\]
Indeed, choose maps $j_X:X\to X^{\vee\vee}$ and
$j_Y:Y\to Y^{\vee\vee}$, and use their tensor product to define $j_{XY}$.
Multiplicativity of the quantum traces defining the squared norm gives
the formula; see \cite[Proposition~4.7.3 and equation~(7.70)]{EGNO}.
For $m_{st}=3$, the words $Y_sY_t$, $Y_tY_s$ and
$Y_sY_tY_s\cong Y_tY_sY_t$ are simple. Hence
\[
 n(Y_s)^2n(Y_t)=n(Y_sY_tY_s)
 =n(Y_tY_sY_t)=n(Y_t)^2n(Y_s),
\]
so $n(Y_s)=n(Y_t)$. The two possible values distinguish the seed
branches, proving the assertion. This argument takes place in the finite
dihedral subcategory generated by $Y_s,Y_t$.

Write $\Czero^+(W)=\Czero(W)$ and let $\Czero^-(W)$ be the category
obtained by applying the field embedding $d_+\mapsto d_-$ to its complete
algebraic structure data. The same construction identities hold after
this embedding. For uniform Yang--Lee markings, choose an automorphism of $\CC$
whose restriction to the field of seed coefficients takes the
Yang--Lee data to the fixed positive data. Twist the candidate,
the reference, and their markings by this same automorphism.
Theorem~\ref{thm:main} gives a complex-linear marked equivalence
between the twisted categories. Twisting both its source and target
back cancels the scalar-field change, so the resulting equivalence
between the original Yang--Lee categories is again complex-linear.
This proves uniform Yang--Lee marked uniqueness.
The positive based-ring dimension in either case is
$h_w\mapsto\ph^{\ell(w)}$; the conjugate spherical dimension is
$d_-^{\ell(w)}$.

\begin{corollary}[Branch choices in finite ADE type]\label{cor:ADE-branches}
Let $W=\prod_{i=1}^r W_i$ be of finite ADE type, with each $W_i$
connected. Consider marked categorifications of $\HH(W)$ whose
one-generator tensor seeds are chosen from $\Fib$ and the Yang--Lee category.
Such a categorification exists precisely when the chosen branch is
constant on each connected component. For each assignment
$\epsilon_i\in\{+,-\}$, it is unique up to marked complex-linear tensor
equivalence and is represented by
\[
 \boxtimes_{i=1}^r\Czero^{\epsilon_i}(W_i).
\]
\end{corollary}
\begin{proof}
Every edge in a connected finite ADE diagram has exponent three, so the
squared-norm argument makes its seed branch constant. Let $\cD_i$ be
the full fusion subcategory on $W_i$. It has the indicated uniform
branch, and the positive-seed theorem or its Galois conjugate identifies
it with $\Czero^{\epsilon_i}(W_i)$.

It remains to identify how the component subcategories combine.
The product basis of $\HH(W)$ first gives an equivalence of linear
categories: the multiplication functor from their Deligne product
to the candidate is an equivalence because products of component simples are simple, distinct, and
exhaust all labels. Objects in different components commute up to
isomorphism, since the cross-component exponents are two. Every
intermediate subcategory has trivial universal grading group: a
homogeneous generator $X$ with $X^2=\one+X$ has degree $g$ satisfying
$g^2=e=g$. These properties supply the hypotheses of the tensor splitting
theorem \cite[Theorem~3.15]{ENOtwisted}. Applied successively, it
upgrades the component decomposition to a tensor equivalence that
is tautological on every component subcategory. Composing with
the marked component equivalences proves uniqueness. Conversely, the
displayed Deligne product realizes every componentwise branch choice.
Different assignments give different squared norms at some labelled
generator, hence inequivalent marked categories. Thus there are $2^r$
marked classes when all branch choices are allowed.
\end{proof}

For general Coxeter systems in \eqref{eq:exponents}, the branch is
constant on each connected component of the graph formed by the
exponent-three edges. An exponent-$\infty$ edge permits different
branches: Yamagami's algebraic free product
$\Fib*\mathrm{YL}$ has alternating simple words and realizes
$\HH(I_2(\infty))$ \cite[Section~3]{Yamagami}.

\subsection{Infinite type and finite-support realization}
For finite $W$, the categories above have $|W|$ simple classes.
For infinite $W$, each object still has finite support and each
edge scalar is a finite integer monomial in the local
parameters. The resulting module autoequivalences have invertible
structure maps on every object, without any infinite summation or
completion. When the actions and markings are unitary, these maps
are unitary by \cref{cor:unitary-comparison}.

The positive-seed classification therefore applies to every system
in \eqref{eq:exponents}, without a finiteness assumption on $W$. The
uniform Yang--Lee uniqueness argument of \cref{subsec:branches}
also applies in infinite type.

\subsection{The common gauge matrix}
Assume $W$ is finite, and put $N=|W|$, $n=|S|$. Let $G$ be the integer matrix whose rows are the exponents in \eqref{eq:gauge-k} and \eqref{eq:gauge-h}, with columns indexed by left and right upward coloured edges. Equivalently,
\[
 G=(\partial_2^{\mix})^{\mathsf t}
\]
in the indicated cellular bases.

\begin{corollary}[Gauge rank and saturation]\label{cor:gauge-rank}
\[
 \rank G=(n-1)N+1.
\]
Every nonzero Smith invariant of $G$ is $1$.
\end{corollary}
\begin{proof}
There are $N$ vertices and $nN$ edges: for each of the $n$ colours there are $N/2$ left edges and $N/2$ right edges. The graph is connected, so its integral cycle group has rank $nN-N+1$. By \cref{thm:pure-elimination} with $k=2$ and the vanishing of $H_1(Y;\ZZ)$,
\[
 \im\partial_2^{\mix}=\im\partial_2=\ker\partial_1.
\]
The latter is primitive in the free edge group, since the quotient is the image of $\partial_1$, a free abelian group. This proves the rank formula and the assertion about Smith invariants, also for the transpose.
\end{proof}

Saturation means that the nonzero part of the matrix of edge rescalings has
no nontrivial Smith factors. Thus its scalar equations can be solved
by integer monomials, with no additional ambiguity from roots of
unity. This is the matrix form of the integral normalization in
\cref{sec:global}.

\section{Discussion}\label{sec:discussion}

\subsection{The role of ADE and Coxeter systems}
The ADE classification obtained here belongs to a broad family of
classification results involving the same diagrams. Finite ADE
diagrams occur as the irreducible simply-laced types in the
Killing--Cartan classification and the classification of finite
Coxeter systems \cite{Bourbaki,Coxeter1934}, and as the finite
connected simple graphs of spectral radius less than two
\cite{Smith,GHJ}. They also appear in quiver representation theory
\cite{Gabriel} and the theory of simple surface singularities
\cite{Slodowy}. The classical McKay correspondence instead produces
affine ADE diagrams from tensoring irreducible representations of
finite subgroups of $\mathrm{SU}(2)$ with the natural two-dimensional
representation \cite{McKay}.

Within tensor category theory, closely related ADE phenomena arise
from subfactors and quantum McKay correspondence. Jones's index
theorem \cite{JonesIndex} underlies the classification of subfactors
of index below four, where the principal graphs are governed by
ADE \cite{GHJ}. In the module-category setting, ADE diagrams encode
the action of the fundamental object on indecomposable module
categories over the $\mathfrak{sl}_2$ fusion categories at positive
integral levels \cite{OstrikModules}. Quantum McKay correspondence relates this
picture to quantum subgroups and conformal field theory
\cite{KirillovOstrik}. Further related classifications concern ADE
modular invariants \cite{CIZ} and local conformal nets with central
charge $c<1$ \cite{KawahigashiLongo}.

The ADE diagrams in our classification record relations among the
generating Fibonacci subcategories, while the simple objects are
indexed by the elements of the associated Coxeter group. Thus the
vertices of our diagrams specify the seeds to be composed, and the
edges encode their pairwise braid relations. For example, the
three vertices of $A_3$ describe three Fibonacci seeds, whereas the
resulting category has $|W(A_3)|=24$ simple objects. This differs from
the role of vertices as simple-object labels in a principal graph
or a module-category action graph.

The Coxeter system is essential throughout the argument. It is
recovered from the based-ring axioms and the Fibonacci relation in
\cref{thm:ring-rigidity}, rather than supplied as an indexing diagram.
Coxeter length, reduced expressions, and parabolic subgroups organize
the simple objects, multiplication, and extension of local braid
maps in \cref{sec:braids,sec:construction}. The geometry of the Davis
complex and its parabolic cells then supplies the global
compatibility needed to normalize the associativity data and prove
marked uniqueness in \cref{sec:geometry,sec:cochain,sec:integral,sec:global}.
Finally, the categorifiability restriction leaves only exponents
$2,3,\infty$, and the simple-object indexing identifies finite
categorical rank with finiteness of $W$. The classification of finite
Coxeter systems then gives precisely the finite ADE diagrams.
The role of Coxeter theory therefore extends from the fusion rules
to the tensor structure and its uniqueness, before the finite-type
criterion is applied.

\subsection{Composites of other seed categories}
The broader composition program concerns arbitrary finite families
$(\mathcal A_i)_{i\in I}$ of prescribed seed tensor categories,
including repeated copies of one seed and mixtures of different
seeds. The data to be classified are a tensor category $\mathcal D$
and fully faithful tensor functors
$\iota_i:\mathcal A_i\to\mathcal D$ whose images generate
$\mathcal D$ under tensor products, finite direct sums, and summands.
An equivalence between two such composites is a tensor equivalence
$F:\mathcal D\to\mathcal D'$ together with monoidal natural
isomorphisms $F\circ\iota_i\cong\iota_i'$ for each $i\in I$.
The tensor structures on the seeds are part of the prescribed data.

A natural next family of seeds is provided by the even parts of
$A_n$ subfactors \cite{GHJ}. For each seed one fixes the subfactor
and the chosen even part, including its tensor structure. The
$A_4$ case recovers the positive Fibonacci seed considered here.
One can then study composites of any finite number of copies of a
fixed even part, as well as composites involving even parts for
different values of $n$. This places the present classification
within a sequence of composition problems with more complicated
seed fusion rules.

Several questions arise already for a pair of seeds: which mixed
fusion rules can occur, which admit categorification, and which
produce finite-rank composites? For an arbitrary finite family,
there is the further problem of deciding when prescribed pairwise
composites can be realized simultaneously. Relations involving
three or more seeds may impose conditions that are invisible in
any single pair. Even when a simultaneous realization exists, one
must determine whether the fixed seed structures and mixed fusion
rules determine its associator up to marked tensor equivalence,
or whether additional parameters are required.

In the Fibonacci case, the quadratic fusion relation forces a
Coxeter system, and its parabolic geometry controls the compatibility
of the local tensor data. For finite families of other seeds, a
natural direction is to investigate the algebraic and geometric
structures underlying their mixed fusion rules and associativity
constraints.

\appendix
\section{The complete local polynomial certificate}\label{app:existence}

This appendix specifies the finite calculation in \cref{lem:local-braid}.
The matrices express the right-module naturality of the proposed
left braid. We fix their path coordinates, list the nonzero polynomial
residuals, and then specialize the seed parameter; the passage to
general parabolic blocks is proved in \cref{prop:global-braid}. For this appendix alone the Coxeter generators are called $a,b$; they correspond to $s,t$ in \cref{sec:A2}, and are unrelated to the scalar entries in \eqref{eq:U}. All formulas are over $\ZZ[z]$ until
$z$ is specialized. The group is
$W=\{e,a,b,ab,ba,o\}$ with $a^2=b^2=e$ and $aba=bab=o$.
Its length function is $0,1,1,2,2,3$ in that order.

\subsection{Path bases and elementary maps}
An operation is either $L_s$ or $R_r$. For an input vertex $w$, its
ordered possible outputs are respectively
\[
\mathcal N_{L_s}(w)=
\begin{cases}(sw),&\ell(sw)>\ell(w),\\(sw,w),&\ell(sw)<\ell(w),\end{cases}
\]
and the analogous list $(wr)$ or $(wr,w)$ for $R_r$.
For an operation list $(O_1,\ldots,O_k)$ in execution order, recursively
enumerate paths $(w_0,w_1,\ldots,w_k)$ by taking
$w_i\in\mathcal N_{O_i}(w_{i-1})$ in the displayed order, exhausting
all continuations of the first choice before the second choice.
Matrix columns are source paths and rows are target paths.

For the interchange of adjacent operations $(R_r,L_s)$ into
$(L_s,R_r)$, an entry is zero unless the starting and ending vertices
coincide. With fixed endpoints, the block is $1$, except in
\eqref{eq:exception}, where it is $U(z)$ and the middle
vertices have order $(sw,w)$. For the braid
$(L_a,L_b,L_a)\to(L_b,L_a,L_b)$, fixed-endpoint blocks are $1$
except for the block specified in
\eqref{eq:source-braid}--\eqref{eq:target-braid}, which is $U(z)$.
These rules determine every matrix used below.

For a right colour $r$, the source operation list for
\eqref{eq:canonical-naturality} is
\[
(R_r,L_a,L_b,L_a),
\]
and the final target list is
\[
(L_b,L_a,L_b,R_r).
\]
For the left side, apply the braid to the last three initial operations
and then move $R_r$ across the resulting three left operations. For
the right side, first move $R_r$ across all three left operations and
then apply the braid. Each move replaces only its indicated subpath
and fixes the remaining vertices. Summing over new middle vertices
is ordinary matrix composition.

\subsection{Matrix sizes and residuals}
The sizes of the square matrices, with all endpoint blocks combined,
are given in Table~\ref{tab:matrix-sizes}. Their squared sizes sum to
$470$.
\begin{table}[htbp]
\centering
\caption{Dimensions in the twelve local module-naturality equations.}
\label{tab:matrix-sizes}
\begin{tabular}{@{}c*{6}{r}@{}}
\toprule
Input $u$ & $e$ & $a$ & $b$ & $ab$ & $ba$ & $o$\\
\midrule
Right colour $a$ & 2&3&4&7&6&11\\
Right colour $b$ & 2&4&3&6&7&11\\
\bottomrule
\end{tabular}
\end{table}

Table~\ref{tab:residuals} lists every entry for which the two sides of
\eqref{eq:canonical-naturality} differ as polynomials. Row and column
indices start at $1$ and refer to the path order just defined.
Every other one of the $460$ entries agrees in $\ZZ[z]$.
Writing $P=z^4+z^2-1$, the final column is the difference of the two
preceding columns divided by $P$.
\begin{table}[htbp]
\centering
\caption{All nonzero formal residuals before imposing $P(z)=0$.}
\label{tab:residuals}
\small
\begin{tabular}{@{}ccrrccc@{}}
\toprule
$u$ & $r$ & Row & Column & Left side & Right side & Residual$/P$\\
\midrule
$ab$ & $a$ & 5 & 3 & $z^4+z^2$ & $1$ & $1$\\
$ab$ & $a$ & 7 & 7 & $z^4+z^2$ & $1$ & $1$\\
$o$ & $a$ & 9 & 11 & $z^5+z^3$ & $z$ & $z$\\
$o$ & $a$ & 11 & 7 & $z^5+z^3$ & $z$ & $z$\\
$o$ & $a$ & 11 & 11 & $-z^6+z^2$ & $z^4$ & $-z^2$\\
$o$ & $b$ & 4 & 3 & $1$ & $z^4+z^2$ & $-1$\\
$o$ & $b$ & 5 & 11 & $z$ & $z^5+z^3$ & $-z$\\
$o$ & $b$ & 10 & 10 & $1$ & $z^4+z^2$ & $-1$\\
$o$ & $b$ & 11 & 4 & $z$ & $z^5+z^3$ & $-z$\\
$o$ & $b$ & 11 & 11 & $z^4$ & $-z^6+z^2$ & $z^2$\\
\bottomrule
\end{tabular}
\end{table}

For example, for $(u,r)=(ab,a)$ the entry in row $5$, column $3$
compares the source path
\[
(ab,o,ba,ba,o)
\]
with the target path
\[
(ab,o,o,ab,o).
\]
The two compositions give $z^4+z^2$ and $1$, respectively.
All remaining entries follow by the same replacement rule; the tables
and the elementary block formulas specify the calculation without
numerical approximation. Finally,
\[
U(z)^{\mathsf t}U(z)=(z^4+z^2)I.
\]
Consequently all required algebraic identities hold at $z=\ph^{-1/2}$. The same specialization makes $U(z)$ unitary by
\eqref{eq:seed-identities}; \cref{prop:unitary-realization} gives the resulting
unitary tensor realization.

\section{Explicit rank-two algebraic certificates}\label{app:rank-two}
This appendix supplies the coefficients used in \cref{prop:rank-two}.
Its three steps mirror that proof: normalize all but one local
parameter by edge rescalings, propagate the existing braid
from its value at the unit, and impose its right-module naturality.
The normalization matrices and every resulting nonzero braid block
are displayed below. The calculations use actual inverse interchange matrices over $\Cx$. To distinguish Coxeter letters from the scalar entries $a,b$ of $U$, we denote the two generators by $s,t$. All computations take place in
\begin{equation}\label{eq:laurent-ring}
R=\ZZ[z,q,q^{-1}]/(z^4+z^2-1),\qquad a=z^2,\quad b=z.
\end{equation}
The positive specialization is $z=\ph^{-1/2}$. The indeterminate $q$ is invertible and otherwise unrestricted.

\subsection{Integral gauge normalization}
Write $k_{s,r}(x)$ or $h_{s,r}(x)$ for the parameter of the two- or four-point orbit of the ordered colour pair $(s,r)$ with shortest vertex $x$. Edge variables $A_s(x),B_r(x)$ are oriented upward. Each matrix entry below records an integer exponent in \eqref{eq:gauge-k}--\eqref{eq:gauge-h}.

For $m_{st}=2$, order $W=\{e,s,t,st\}$ and the six parameters as
\begin{equation}\label{eq:rank2-phase-list}
\bigl(k_{s,s}(e),k_{s,s}(t),h_{s,t}(e),h_{t,s}(e),
      k_{t,t}(e),k_{t,t}(s)\bigr).
\end{equation}
On the first five rows and columns
\[
\bigl(A_s(e),A_s(t),B_s(e),A_t(s),A_t(e)\bigr)
\]
the gauge matrix is
\begin{equation}\label{eq:M2}
M_2=\begin{pmatrix}
-1&0&1&0&0\\
0&-1&0&0&0\\
-1&1&0&0&0\\
0&0&1&1&-1\\
0&0&0&0&-1
\end{pmatrix},\qquad \det M_2=1.
\end{equation}
Hence the first five parameters can be made one. Denote the remaining value of $k_{t,t}(s)$ by $q$.

For $m_{st}=3$, let $W=\{e,s,t,st,ts,o\}$, with $o=sts=tst$, and order the eight parameters as
\begin{equation}\label{eq:rank3-phase-list}
\begin{split}
(&k_{s,s}(e),h_{s,s}(t),h_{s,t}(e),k_{s,t}(ts),\\
 &h_{t,s}(e),k_{t,s}(st),k_{t,t}(e),h_{t,t}(s)).
\end{split}
\end{equation}
For the first seven rows and columns
\[
\bigl(A_s(e),A_s(t),A_s(ts),B_s(e),A_t(s),A_t(st),A_t(e)\bigr)
\]
the matrix is
\begin{equation}\label{eq:M3}
M_3=\begin{pmatrix}
-1&0&0&1&0&0&0\\
0&-1&1&0&0&0&0\\
-1&1&0&0&0&0&0\\
0&0&-1&0&0&0&0\\
0&0&0&1&1&0&-1\\
0&0&0&0&0&-1&0\\
0&0&0&0&0&0&-1
\end{pmatrix},\qquad \det M_3=-1.
\end{equation}
The first seven parameters can therefore be set to one; put $q=h_{t,t}(s)$ afterward. In either case the inverse of $M_i$ has integer entries, so the required gauges are integer monomials in the original nonzero scalars. These integer-monomial gauges fix the parameters exactly, including their finite-order parts.

\subsection{Path ordering and the forced braid}
For $L_s$ at $w$, list outputs as $(sw)$ on an ascent and $(sw,w)$ on a descent; for $R_r$ replace $sw$ by $wr$. Enumerate paths recursively in this order, exhausting the continuations of the first choice before those of the second. Rows index target paths and columns source paths. Lists of operations are written in execution order, in contrast to composition of functors.

The source and target left-operation lists are
\[
\begin{array}{c|cc}
m&\text{source}&\text{target}\\\hline
2&(L_s,L_t)&(L_t,L_s)\\
3&(L_s,L_t,L_s)&(L_t,L_s,L_t).
\end{array}
\]
Thus for $m=2$ the source functor is $L_tL_s$. Normalize $T_{\eobj{e}}=1$ and propagate its components by \eqref{eq:forced-braid}. The right words used for the inputs $s,t,st$ at $m=2$ are $s,t,ts$; those for $s,t,st,ts,o$ at $m=3$ are $s,t,st,ts,sts$. They are execution-order lists and give the required input from $e$.

The resulting nonzero blocks are displayed below. An entry $(v,c)$ denotes the coefficient of the block with final vertex $v$; all undisplayed entries are zero. Except for the displayed two-dimensional block, the relevant source and target path spaces are one-dimensional.

\begin{table}[ht]
\centering
\caption{Forced braid blocks for $m=2$.}\label{tab:T2}
\begin{tabular}{cl}
\toprule
Input&Nonzero blocks $(\text{final vertex},\text{coefficient})$\\\midrule
$e$&$(st,1)$\\
$s$&$(t,1),\ (st,1)$\\
$t$&$(s,q^2),\ (st,q)$\\
$st$&$(e,q^2),\ (s,q^2),\ (t,q),\ (st,q)$\\
\bottomrule
\end{tabular}
\end{table}

\begin{table}[ht]
\centering
\caption{Forced braid blocks for $m=3$.}\label{tab:T3}
\begin{tabular}{cl}
\toprule
Input&Nonzero blocks $(\text{final vertex},\text{coefficient})$\\\midrule
$e$&$(o,1)$\\
$s$&$(st,1),\ (o,1)$\\
$t$&$(ts,q^2),\ (o,q)$\\
$st$&$(s,q^{-2}),\ (st,q^{-1}),\ (ts,1),\ (o,1)$\\
$ts$&$(t,q^2),\ (st,q),\ (ts,q^2),\ (o,q)$\\
$o$&$(e,q^{-2}),\ (s,q^{-2}),\ (t,1),\ (st,1),\ (ts,1)$\\
&and the $(o,o)$ block $\displaystyle\begin{pmatrix}a/q&b\\b/q&-a\end{pmatrix}$\\[8pt]
\bottomrule
\end{tabular}
\end{table}

For the last block the source paths are
\[
(o,ts,ts,o),\qquad(o,o,o,o),
\]
and the target paths are
\[
(o,st,st,o),\qquad(o,o,o,o),
\]
in those orders. Move the right word past the left word using the elementary mixed interchanges, identify the normalized value at $e$, and reverse the target propagation using \eqref{eq:Uk-inverse}. This computes the tables directly. Reduction of coefficients uses only $a=b^2$ and $a^2+b^2=1$; it never uses $\overline q=q^{-1}$.

\subsection{The residual module-naturality equations}
For \eqref{eq:braid-natural}, start with the right operation followed by the source left list. One composite applies $T$ first and then moves the right operation past the target list; the other moves it past the source list first and then applies $T$. The following entries of these two matrices give the necessary equations. Row and column numbers are one-based in the preceding recursive order.

\begin{table}[ht]
\centering
\caption{Two coefficients forcing the residual parameter.}\label{tab:rank-two-witness}
\begin{tabular}{cccccc}
\toprule
$m$&Input&Right colour&Row, column&Left side&Right side\\\midrule
$2$&$s$&$t$&$4,4$&$q^2$&$q$\\
$3$&$ts$&$t$&$5,5$&$1$&$q$\\
\bottomrule
\end{tabular}
\end{table}

For $m=2$ the source and target paths are both $(s,st,st,st)$, with operation lists $(R_t,L_s,L_t)$ and $(L_t,L_s,R_t)$. The only nonstandard elementary orbit is the two-point $(t,t)$ orbit with shortest vertex $s$.

For $m=3$, the source path is
\[
(ts,o,o,st,st)\quad\text{for }(R_t,L_s,L_t,L_s),
\]
and the target path is
\[
(ts,ts,o,st,st)\quad\text{for }(L_t,L_s,L_t,R_t).
\]
The only nonstandard elementary orbit is the four-point $(t,t)$ orbit with shortest vertex $s$. The left coefficient in the second row simplifies to one using $a^2+b^2=1$. Consequently
\[
q(q-1)=0\quad(m=2),\qquad q-1=0\quad(m=3).
\]
Because $q$ is invertible, both equations force $q=1$. Both equations follow from coefficients of the candidate's module-naturality identity.

\subsection{Reproducibility of the finite calculations}
The script \path{Fibonacci_Coxeter_Construction_Uniqueness_Checks.py}
uses the Python standard library to enumerate Coxeter elements and
ordered paths and compose elementary matrices exactly.
For \cref{app:existence}, it reconstructs the twelve matrices over
$\ZZ[z]$, checks their sizes, and reproduces the complete nonzero
polynomial residual list. For \cref{app:rank-two}, it works over
$\ZZ[z,q,q^{-1}]/(z^4+z^2-1)$, verifies inverse interchanges,
propagates the normalized braid, and obtains $q^2=q$ and $1=q$.
It also checks every generator module-naturality entry at $q=1$
and both integer gauge determinants.

The scripts \path{A3_Word_Checks.py} and
\path{A3_Associator_Counts.py} check, respectively, the permutation
labels and inversion counts in \eqref{eq:A3-five}--\eqref{eq:A3-prefixes}
and the associator-table counts in \cref{subsec:A3}.
All three scripts and their JSON outputs are included in \path{anc/};
the matrix script writes its output beside the script.

These computations verify the displayed finite identities.
Their extension to global module naturality is proved in
\cref{prop:global-braid}; the geometric compatibility and integral
normalization are proved in \cref{sec:geometry,sec:cochain,sec:integral,sec:global}.

\section*{Acknowledgments}
Both authors acknowledge support from the Beijing Natural Science
Foundation Key Program (Grant~Z220002) and the Beijing Natural Science
Foundation (Grant~Z221100002722017).

\par\medskip
\noindent\textbf{AI-use disclosure.}
The AI tool was used as a research assistant.
The authors retain full responsibility for the final
version of the manuscript.

\end{document}

%% file: companion_refs.tex
\newcommand{\PointedArithmeticTheorem}{4.11}
\newcommand{\PointedRigidityTheorem}{6.1}